\documentclass[10pt,a4paper,oneside]{article}
\usepackage[margin=2.5cm,includefoot,footskip=30pt,headheight=13.6pt]{geometry}
\usepackage[utf8]{inputenc}
\usepackage{epsfig,graphics,graphicx,amsmath,amssymb,amsfonts}
\usepackage{epstopdf}
\usepackage{subcaption}
\usepackage{color}
\usepackage{url}
\usepackage{listings}
\usepackage{algorithm}
\usepackage{algpseudocode}
\usepackage{float}
\usepackage{overpic}

\newcommand{\R}{\mathbb{R}}

\newcommand{\D}{\mathcal{D}}

\author{Federico Pichi, Annalisa Quaini, Gianluigi Rozza}
\title{A Reduced Order Modeling technique to study bifurcating phenomena: application to the Gross-Pitaevskii equation}

\newcommand{\br}{\mathbf{r}}
\newcommand{\bx}{\mathbf{x}}
\newcommand{\mmu}{\boldsymbol{\mu}}

\def\cl {\nonumber \\}
\def\el {\nonumber}

\begin{document}
\maketitle

\noindent {\bf Abstract.} We propose a computationally efficient framework to treat nonlinear
partial differential equations having bifurcating solutions as 
one or more physical control parameters are varied.
Our focus is on steady bifurcations.
Plotting a bifurcation diagram entails computing multiple solutions of a parametrized, nonlinear problem,
which can be extremely expensive in terms of computational time.
In order to reduce these demanding computational costs, our approach 
combines a continuation technique and Newton's method with
a Reduced Order Modeling (ROM) technique, suitably supplemented with
a hyper-reduction method.
To demonstrate the effectiveness of our ROM approach, we trace the steady solution branches
of a nonlinear Schr\"{o}dinger equation, called Gross--Pitaevskii equation, as
one or two physical parameters are varied. In the two parameter study, we show that our approach  is
60 times faster in constructing a bifurcation diagram than a standard Full Order Method.

\section{Introduction}\label{sec:intro}

We consider the problem of finding a solution $X \in V$ such that
\begin{align}\label{eq:gen_problem}
G(X(\mmu); \mmu) = 0,
\end{align}
where $\mmu \in \D $ denotes a point in a parameter domain $\D \subset \mathbb{R}^M$, $V$ a
given functional space, $G(\cdot,\mmu)$ a given nonlinear functional in $X$. 
We view problem \eqref{eq:gen_problem} as the strong form of a nonlinear partial differential equation (PDE) or a system of such equations
in which $M$ parameters appear.
We are interested in situations in which the solution $X$ of \eqref{eq:gen_problem} differs in character 
for parameter vectors $\mmu$ in different subregions of the parameter domain $\D$. 
Such situations occur if $X$ undergoes bifurcations as $\mmu$ changes from one subregion to another.

In general, one cannot solve \eqref{eq:gen_problem} for $X$ so that one instead seeks an approximation of
the solution in a $N_h$-dimensional subspace $V_h \subset V$. 
Such approximation is usually obtained using a so-called Full Order Method (FOM), like for example
the Finite Element Method, which is often expensive, especially if multiple solutions are needed. 
For this reason, one is interested in finding surrogates methods that are much less costly 
so that obtaining approximations to the exact solution $X$  of \eqref{eq:gen_problem} 
for many choices of the parameter vector $\mmu\in \D$ 
becomes feasible. Such surrogates, called Reduced Order Modeling (ROM) techniques,
are constructed using a ``few'' solutions computed with the FOM. 

A practical way to realize a ROM is to organize the computation in two steps:
\begin{itemize}
\item[-] An {\bf offline phase}: approximation solutions corresponding 
to selected representative parameters values/system configurations are computed 
with a FOM and stored, together with other information about the parametrized problem. 
This is a computationally expensive step usually performed on high performance computing facilities.
\item[-] An {\bf online phase}: the information obtained during the offline phase is used to compute the solution for a newly specified value of the parameters in a short amount of time (ideally in real time), even on a relatively low power device such as a laptop or a smartphone.
\end{itemize}
These split computational procedures are built in such a way that new parameter dependent quantities are easily and quickly computed online, while representative basis functions for selected parameter values and more demanding quantities are pre-computed offline.
Among all the possible ROMs, we choose the Reduced Basis (RB) technique
\cite{hesthaven2015certified,patera07:book,quarteroni2015reduced}.

Recent developments of ROM techniques have focused on the reduction of
computational time for a wide range of differential problems \cite{RozzaEncyclopedia,QuarteroniRozza2013}, 
while maintaining a prescribed tolerance on
error bounds \cite{Rozza:ARCME,hesthaven2015certified, quarteroni2015reduced}.
ROMs in the setting of bifurcating solutions are considered in the early papers 
\cite{NOOR1982955,Noor:1994,Noor:1983,NOOR198367} for buckling bifurcations in solid mechanics. 
More recently, in \cite{Terragni:2012} it is shown that a POD approach allows for considerable computational time savings for the analysis of bifurcations in some nonlinear dissipative systems. 
RB methods have been used to study symmetry breaking bifurcations
\cite{Maday:RB2,PLA2015162} and Hopf bifurcations \cite{PR15} for
natural convection problems. A RB method for symmetry breaking bifurcations in contraction-expansion channels
has been proposed in  \cite{PITTON2017534}. In \cite{HESS2019379}, steady bifurcations for both a natural convection
problem and contraction-expansion channels are investigated with a localized ROM approach.
Yano and Patera \cite{Yano:2013} introduced a RB
method for the stability of flows under perturbations in the forcing term or in the boundary conditions,
which is based on a space-time framework that allows for particularly sharp error estimates. Furthermore, in~\cite{Yano:2013} 
it is shown how a space-time inf-sup constant approaches zero as the computed solutions get close to a bifurcating value. 
A recent work on ROMs for bifurcating solutions in structural mechanics is \cite{pichirozza}.
Finally, we would like to mention that machine learning techniques based on sparse optimization 
have been applied to detect bifurcating
branches of solutions in \cite{BTBK14,KGBNB17} for a two-dimensional laterally heated cavity and Ginzburg-Landau model, respectively. 

The methodology we propose to plot a bifurcation diagram makes use of three building blocks: 
i)~a continuation technique to properly follow each solution branch, ii)~Newton's method to 
deal with the nonlinearity of problem \eqref{eq:gen_problem},
 and iii)~a RB method to efficiently solve the linearized problem obtained from Newton's method. 
The novelty of this approach relies on the combination of well-known and assessed 
methods to obtain a versatile and global approach to any kind of bifurcation problems modeled 
by parametric PDEs. We present technical insights to build a ROM approach capable of 
approximating efficiently and accurately the whole bifurcation diagram with the use of a unique reduced basis.
As a concrete setting to illustrate our methodology, we consider the Gross--Pitaevskii equation.

Often referred to as a nonlinear Schr\"{o}dinger equation, 
the Gross--Pitaevskii equation models certain classes of Bose-Einstein condensates
(BECs), a special state of matter formed by identical  bosons at ultra-low temperatures.
It is well known that the solutions of the Gross--Pitaevskii equation with a parabolic trap
in two dimensions exhibit a rich bifurcating behavior
\cite{Middelkamp_et_al2010,Middelkamp_et_al2011,Charalampidis_et_al2018}, 
which includes symmetry-breaking bifurcations
and vortex-bearing states when a (sufficiently strong) rotational angular momentum term is added \cite{GarciaAzpeitia2017}.
The bifurcating behavior becomes even richer for the 
two-dimensional coupled Gross-Pitaevskii equations \cite{Charalampidis_et_al2019}.
However, for simplicity we stick to the simple Gross-Pitaevskii equation and present a one parameter study 
(the chemical potential being the only varying parameter) 
and a two parameter case (varying chemical potential and the normalized trap strength).
For both cases, we show that our approach is able to capture the first six solution branches with high accuracy.
We stress the need to supplement the RB method with a hyper-reduction technique \cite{barrault04:_empir_inter_method,Chaturantabut2010}
to enable significant computational time savings with respect to a standard Full Order Method.
In particular, we show that if the RB method at the above point iii) is combined with the 
Discrete Empirical Interpolation Method \cite{Chaturantabut2010}, our approach is up to
60 times faster in constructing a bifurcation diagram than a FOM, making it an ideal tool
to study the complex solution behavior of the Gross--Pitaevskii equation
and other nonlinear PDE problems.

The work outline is as follows. In Sec.~\ref{sec:num_approx}, we present the 
building blocks of our approach to reduce the computational time
requested by the construction of the bifurcation diagrams. In Sec.~\ref{sec:gp_eq}, we apply
such approach to the Gross--Pitaevskii equation.
Numerical results pertaining to the validation of the FOM, reconstruction of bifurcation 
diagrams with our ROM approach with and without hyper-reduction technique
are reported in Sec.~\ref{sec:num_res}. Conclusions are drawn in Sec.~\ref{sec:concl}.

\section{Numerical approximation of a problem with bifurcations}\label{sec:num_approx}

The nonlinearity in problem \eqref{eq:gen_problem} 
can produce a loss of uniqueness for the solution, with multiple solutions
branching from a known solution at a bifurcation point.
Our aim is to study numerically the associated bifurcation diagrams
with contained computational costs. We restrict our attention
to steady bifurcations.

The standard assumption for the map $G$ in eq~\eqref{eq:gen_problem}
 is the continuous differentiability with respect to $X$ and $\mmu$. Let $(\bar{X}, \bar{\mmu}) \in V \times \mathcal{D}$ be the known solution, i.e.~$G(\bar{X}, \bar{\mmu}) = 0$. Let us denote by $D_XG(Z,\mmu): V \rightarrow V'$ and $D_{\mmu}G(Z,\mmu): \mathcal{D} \rightarrow V'$ the partial derivatives of $G$ on a generic point $(Z, \mmu) \in V \times \mathcal{D}$. 
 A strong assumption usually found in the literature in order to have a local branch of non-singular solutions 
is that $D_XG(\bar{X}, \bar{\mmu}): V \rightarrow V'$ is bijective. Of course, this is not our case: 
we deal with bifurcation points, which are singularities for the system. Moreover, we do not require $G$
to be affine in $\mmu$, because of the nonlinearity of the problem. As we will see later, this forces us 
to implement a hyper-reduction technique, such as EIM/DEIM  \cite{barrault04:_empir_inter_method,Chaturantabut2010}, to recover efficiency.

\subsection{Proposed approach for the branch reconstruction}

For simplicity, we present our approach for a scalar parameter $\mu$ (i.e., $M = 1$)
in parameter domain $\D = [\mu_i, \mu_f ]$, 
although it can easily be extended to a multi-parameter setting (i.e., $M > 1$).
Our algorithm is based on three building blocks:
\begin{enumerate}
\item A \emph{continuation technique}: it permits to reconstruct properly the bifurcation diagram by
following each branch and providing a suitable initial guess for the nonlinear iterations (next building block). We use 
a slight modification of the continuation method in \cite{allogwer}, which consists in a 
for loop in $\bar{\D} = [\mu_1, \dots, \mu_K] \subset \D$, a discrete version of the parameter set with cardinality $K$. 
At each new cycle, we check if a bifurcation occurs, using a threshold $\epsilon_{BIF}$ that controls the norm of the solution. 
If a bifurcation does occur, we set that bifurcated solution as the initial guess for the next cycle 
in order to capture the post-bifurcation behavior.
\item \emph{Newton's method}: each point in the bifurcation diagram requires the solution of nonlinear PDE problem for the corresponding value $\mu$. To deal with the nonlinearity, we use 
the Newton-Kantorovich method \cite{ciarlet2013linear}.
\item A \emph{Galerkin Finite Element discretization}: at each step of the Newton iteration,
we have to approximate the solution of a new linear weak formulation. 
For the space discretization, we use the Galerkin-Finite Element method
as discussed in Sec.~\ref{sec:FOM}. 
This choice is made also in view of the numerical extension towards the model order reduction. 
\end{enumerate}
Algorithm \ref{alg:01} summarizes the three steps for generic nonlinear parametric problem \eqref{eq:gen_problem}. 

\begin{algorithm}
\caption{A pseudo-code for the on-line reconstruction of a branch}\label{alg:01}
\begin{algorithmic}[1]
\For{$j = 1 : K $}\Comment{Continuation loop on $\bar{D}$}
	\If{$||X_{j-1}||_V < \epsilon_{BIF}$} \Comment{Select initial guess}
		\State{$X_j^{(0)} = X_{guess}$} \Comment{Chosen guess}
	\Else
		\State{$X_j^{(0)} = X_{j-1}$} \Comment{Continuation guess}
	\EndIf
	\While{$||\delta X||_{V} > \epsilon$}\Comment{Newton's method}
		\State{$D_XG(X^{(i)}(\mu); \mu)\delta X = G(X^{(i)}(\mu); \mu)$}\Comment{Galerkin FE method}
		\State{$X^{(i+1)}(\mu) = X^{(i)}(\mu) - \delta X$}
	\EndWhile
	\State{$X_{j} = X_{sol}$}
\EndFor
\end{algorithmic}
\end{algorithm}

Let us clarify lines 2-6 in Algorithm \ref{alg:01}.
 In order to decrease the computational cost also during a generic online phase (i.e, also in the case with no hyper-reduction), 
we avoid applying the continuation method when the bifurcated solutions do not exist yet. To do this, we control the norm of the solution 
and we keep providing the pre-bifurcation guess, which is different for each of the branches we want to approximate, until the solver converges 
to a bifurcated solution. At this point, the continuation is enabled and the computed solution becomes the new guess.

We have applied similar strategies for the numerical study of bifurcations
arising in different contexts, ranging from fluid dynamics \cite{PR15,PITTON2017534,HESS2019379,Hess2019}
to structural mechanics \cite{pichirozza}.

\subsection{Galerkin Finite Element method for a (generic) nonlinear problem}\label{sec:FOM}

In this section, we introduce some standard notion for the discretization 
of generic problem \eqref{eq:gen_problem} with the Galerkin Finite Element (FE) method.

Let $V_h$ be a family of finite dimensional spaces, 
such that $V_h \subset V$. Let $N_h = \text{dim}(V_h)$.
We first cast problem \eqref{eq:gen_problem} in weak form and then for a given parameter $\mu \in \D$ 
seek $X_h(\mu) \in V_h$ that satisfies 
\begin{equation}
\label{eq:weakgal}
\langle G(X_h(\mu); \mu), Y_h \rangle \doteq g(X_h(\mu),Y_h ; \mu) = 0\ , \quad \forall\ Y_h \in V_h. 
\end{equation}

To treat the nonlinearity in $G$, we apply
the Newton-Kantorovich method \cite{ciarlet2013linear,quarteroni2015reduced}, which reads as follows.
Choose initial guess $X^0_h(\mu) \in V_h$. Then,  for every $k = 0, 1, \dots$:
\begin{itemize}
\item[-] Step 1: Seek the variation $\delta X_h \in V_h$ such that
\begin{equation}
\label{eq:weaknewton}
dg[X_h^k(\mu)](\delta X_h, Y_h; \mu) =  g(X_h^k(\mu),Y_h ; \mu) \ , \quad \forall\ Y_h \in V_h, 
\end{equation}
\item[-] Step 2: Update the solution
$$
X_h^{k+1}(\mu) = X_h^{k}(\mu) - \delta X_h
$$
\end{itemize}
Step 1 and 2 are repeated until the $L^2$-norm of 
the residual falls below a prescribed tolerance $\epsilon$.

We denote with $\{E^j\}_{j=1}^{N_h}$ a basis for $V_h$. 
Newton's method combined with the Galerkin finite element method and applied to 
problem \eqref{eq:weakgal} reads: find $\delta \vec{X}_h \in \R^{N_h}$ such that 
\begin{equation}
\label{linearnewtgal}
\mathbb{J}(\vec{X}_h^k(\mu); \mu) \delta \vec{X}_h = G_h(\vec{X}_h^k(\mu); \mu) \ ,
\end{equation}
where the Jacobian matrix in $\R^{N_h \times N_h}$ is defined as 
\begin{equation}
\label{jacobiandef}
\mathbb{J}(\vec{X}_h^k(\mu); \mu))_{ij} = dg[X_h^k(\mu)](E^j, E^i; \mu) , \quad \text{for all} \quad   i,j = 1, \dots, N_h \, .
\end{equation}

\subsection{The Reduced Basis method}\label{sec:ROM}

Once projected onto a suitable Finite Element space, the parametrized discrete problem
derived in the previous section leads to a very large nonlinear system which has to be solved 
for every parameter $\mu \in \mathcal{D}$. The solution computed with a FE method
represents the so-called \textit{high fidelity approximation}, which is computationally expensive. 
To reduce the computational cost without compromising the accuracy, 
we choose to use a ROM technique called \textit{Reduced Basis (RB) method} \cite{hesthaven2015certified,patera07:book,quarteroni2015reduced}. 
Roughly speaking, this method consists in a projection of the high fidelity problem on a subspace of smaller dimension, constructed with some properly chosen basis functions.

RB methods use the offline-online paradigm introduced in Sec.~\ref{sec:intro}.
In the \textit{offline phase}, we explore the parameter space $\D$ in order to construct a basis 
for the low dimensional manifold, which efficiently approximates the high fidelity space and
where the parametrized solutions lie.
This entails solving $N_{train}$ times the Galerkin high-fidelity problem associated to
$N_{train}$ values of $\mu$ in $\D$. For the numerical results in Sec.~\ref{sec:num_res},
we chose $\{\mu^n\}_{n=1}^{N_{train}}$ to be an ordered sampling of the interval $\D$. 
In the \textit{online phase} the solution is computed through the projection on the low dimensional manifold
in an efficient and reliable way for every $\mu \in \mathcal{D}$ we are interested in. 
The reduced computational cost comes from avoiding to project on the large Finite Element space.
To be precise, we want to construct the reduced problem through 
the projection on a subspace $V_N \subset V_h$  spanned by a 
collection of the \textit{snapshots}, i.e.~solutions of the full order problem 
for selected values of parameter $\mu$, obtained by, e.g., \textit{Proper orthogonal decomposition} (POD) 
or \textit{Greedy} techniques \cite{hesthaven2015certified, patera07:book, quarteroni2015reduced}.  

Let us provide some more details of the online phase for generic problem \eqref{eq:gen_problem}.
For a given $\mu \in \D$, we seek $X_N(\mu) \in V_N$ that satisfies
\begin{equation}
\label{eq:weakred}
g(X_N(\mu),Y_N ; \mu) = 0, \quad \forall\ Y_N \in V_N, 
\end{equation}
where $g(\cdot, \cdot; \mu)$ is defined in \eqref{eq:weakgal}. Just like for the FOM in Sec.~\ref{sec:FOM}, 
we apply the Newton-Kantorovich method. We
choose an initial guess $X^0_N(\mu) \in V_N$ and the for every $k = 0, 1, \dots$:
\begin{itemize}
\item[-] Step 1: Find the variation $\delta X_N \in V_N$ such that
\begin{equation}
\label{eq:weakrednewton}
dg[X_N^k(\mu)](\delta X_N, Y_N; \mu) =  g(X_N^k(\mu),Y_N ; \mu), \quad \forall\ Y_N \in V_N. 
\end{equation}
\item[-] Step 2: Update the solution 
$$
X_N^{k+1}(\mu) = X_N^{k}(\mu) - \delta X_N.
$$
\end{itemize}
 Step 1 and 2 are repeated until the $L^2$-norm of 
the residual falls below a prescribed tolerance $\epsilon$.

Let $\{\Sigma^m\}_{m=1}^N$ be an orthonormal basis (with respect to the inner product  defined on the space $V_h$) for $V_N$, obtained through POD sampling and the Gram-Schmidt procedure during the offline phase. We remark that this basis will be optimal in the $\ell^2$ sense, minimizing over all possible N-dimensional orthonormal bases $W_N$ the errors between the snapshots and their projection through $W_N$.

Then, $V_N = \text{span}\{\Sigma^1, \dots, \Sigma^N\}$ and we can write every $X_N(\mu) \in V_N$ as
\begin{equation}
\label{eq:solreddecomp}
X_N(\mu) = \sum_{m=1}^{N} X_N^{(m)}(\mu)\Sigma^m. 
\end{equation}
We denote with $\vec{X}_N(\mu) = \{X_N^{(m)}(\mu)\}_{m=1}^{N} \in \R^N$ the reduced solution vector.

By plugging \eqref{eq:solreddecomp} into \eqref{eq:weakred} and choosing
$Y_N = \Sigma^n \in V_N$,  for $1 \leq n \leq N$, we obtain the following algebraic system
\begin{equation}
\label{eq:linearredgalerkin}
g\left(\sum_{m=1}^{N} X_N^{(m)}(\mu)\Sigma^m, \Sigma^n ; \mu\right) = 0 \ , \quad n= 1, \dots, N \ .  
\end{equation}
Let 
\begin{equation*}
(G_N(\vec{X}_N(\mu); \mu))_n = g\left(\sum_{m=1}^{N} X_N^{(m)}(\mu)\Sigma^m, \Sigma^n ; \mu\right),
\end{equation*}
be the \emph{residual reduced vector}. We denote with $\mathbb{V}$ the $N_h \times N$ transformation matrix whose elements 
\begin{equation}\label{eq:matrix_V}
(\mathbb{V})_{jm} = \Sigma^m_{(j)}
\end{equation}
are the nodal evaluation of the \textit{m}-th basis function at the \textit{j}-th node.
With this new notation, we can rewrite problem \eqref{eq:linearredgalerkin} as 
\begin{equation*}
\mathbb{V}^TG_N(\mathbb{V}\vec{X}_N(\mu); \mu) = 0. 
\end{equation*}

Finally, we combine Newton's method and the RB technique. At every iteration $k$ of Newton's method
the problem that has to be solved reads as follows: 
find $\delta \vec{X}_N \in \R^{N}$ such that 
\begin{equation}
\label{linearnewtred}
\mathbb{J}_N(\vec{X}_N^k(\mu); \mu) \delta \vec{X}_N = G_N(\vec{X}_N^k(\mu); \mu),
\end{equation}
where $\mathbb{J}_N$ is the $\R^{N \times N}$ reduced Jacobian matrix
\begin{equation*}
\mathbb{J}_N(\vec{X}_N^k(\mu); \mu) = \mathbb{V}^T \mathbb{J}(\mathbb{V}\vec{X}_N^k(\mu); \mu)\mathbb{V}.
\end{equation*}

We remark that eq.~\eqref{linearnewtred} involves the degrees of freedom of the high fidelity problem. 
Because of this, the repeated assembly of the Jacobian compromises the efficiency of the reduced order 
method during the online phase. As we will see later, this issue can be overcome by adopting an affine recovery technique, which 
allows a consistent speed-up of the method by interpolating the nonlinear part of the variational form.


\section{Application to the Gross--Pitaevskii equation}\label{sec:gp_eq}
The Gross--Pitaevskii equation models certain classes of Bose-Einstein condensates. 
A BEC is a special state of matter formed by 
an unlimited number of bosons that ``condense'' into the same energy state
at low temperatures.
A BEC is formed by cooling a gas of extremely low density, 
about one-hundred-thousandth the density of normal air, to ultra-low temperatures
(close to absolute zero).

A quantum system is the environment to be studied in terms of 
wave-particle duality (i.e., all particles exhibit a wave nature and viceversa) and
it involves the wave-function and its constituents, such as the momentum and wavelength.
The Gross--Pitaevskii equation describes the ground state of a quantum system of identical bosons
using two simplifications: the Hartree--Fock approximation and the pseudopotential interaction model. 
In the Hartree--Fock approximation, the total wave-function 
$\Phi_{tot}$ of a system of $N$ bosons is taken as a product of single-particle functions $\Phi$:
\begin{equation}
\Phi_{tot}(\br_1, \br_2, \dots, \br_N) = \prod_{i = 1}^N \Phi(\br_i), \el
\end{equation}
where $\br_i$ is the coordinate of the $i$-th boson.
If the single-particle wave-function satisfies the Gross--Pitaevskii equation,
the total wave-function minimizes the expectation value (i.e., the probabilistic expected value of the result 
of an experiment) of the pseudopotential model Hamiltonian under normalization condition:
\begin{equation}\label{eq:n_bosons}
N = \int_D \rho ~ d\boldsymbol{r}, \quad \rho = |\Phi|^2,
\end{equation}
where $D \subseteq \mathbb{R}^2$ is the domain under consideration
and $\rho$ is interpreted as the particle density.
The Gross--Pitaevskii equation reads: Find the single-particle wave-function 
$\Phi (\br, t) : \overline{D \times \mathbb{R}^+} \rightarrow \mathbb{C}$ such that
\begin{equation}\label{eq:Schrodinger}
i \partial_t \Phi = \ - \frac{1}{2} \Delta \Phi +|\Phi|^2\Phi + W(r)\Phi \quad \text{in}~D,
\end{equation}
where $i$ is the imaginary unit, $r = | \br | = \sqrt{x^2 + y^2}$ is the radial coordinate,
and $W(r) = \frac{1}{2}\Omega^2 r^2$ is the external potential, with $\Omega$ being the 
normalized trap strength, i.e.~the ratio of trappings along and transverse to the plane.
In this paper, we set $\Omega = 0.2$ unless specified otherwise.
Notice that we consider a single well potential.
Eq.~\eqref{eq:Schrodinger} is similar in form to the Ginzburg--Landau equation 
and is sometimes referred to as a nonlinear Schr\"{o}dinger equation.
Obviously, eq.~\eqref{eq:Schrodinger} needs to be supplemented with 
suitable boundary conditions. 

The construction of the steady solution is based on the ansatz: 
\begin{equation}\label{eq:ansatz}
\Phi(\br,t) = \phi(r)\exp(-i\mu t), \quad \phi(r): \overline{D} \rightarrow \mathbb{C},
\end{equation}
where $\mu$ is the chemical potential,
which has to satisfy $\mu \geq \Omega$.
By plugging \eqref{eq:ansatz} into \eqref{eq:Schrodinger}, 
we obtain nonlinear problem 
\begin{equation}\label{eq:N_Schrodinger}
G(\phi; \mu) \doteq \ - \frac{1}{2} \Delta \phi +|\phi|^2\phi + W(r)\phi - \mu \phi = 0.
\end{equation}

It is well known that the solutions of the one-dimensional version of problem \eqref{eq:N_Schrodinger}
exhibit a bifurcating behavior \cite{KIVSHAR2001225,Kevrekidis_et_al2005,Alfimov2007,Coles_et_al2010}, 
which is not particularly rich  though. 
The bifurcations occurring in the two-dimensional problem \eqref{eq:N_Schrodinger}
are far more interesting \cite{Middelkamp_et_al2010,Middelkamp_et_al2011,Charalampidis_et_al2018}. Indeed, 
several secondary bifurcations appear, which include symmetry-breaking bifurcations
and vortex-bearing states \cite{GarciaAzpeitia2017}.
The bifurcation diagram plots the number of bosons $N$ in the BEC \eqref{eq:n_bosons}
as a function of the chemical potential $\mu$.
When $N \to 0$, the nonlinearity of the problem becomes irrelevant and 
the states bifurcate from the respective linear limit.
Starting from this low-density context, we are interested in exploring the solution 
modes for greater values of $\mu$, which make the problem strongly nonlinear.
Since an arbitrary potential can be approximated as a harmonic potential 
at the vicinity of a stable equilibrium point, when $N \to 0$
we can decompose the linear eigenfunction $\phi_{m,n}$ in Cartesian form as being proportional to 
\begin{equation}
\label{herm_guess}
|m, n\rangle \ :=  \ \phi_{m,n}  \sim H_m(\sqrt{\Omega}x)H_n(\sqrt{\Omega}y)e^{-\frac{r^2}{2}\Omega},
\end{equation}
where $H_j$ is the Hermite polynomial with $j$ being the associated quantum number 
of the harmonic oscillator. 
The critical value of the eigenvalue corresponding to linear eigenfunction $\phi_{m,n}$
is $\mu_{crit} = E_{m,n} :=  (m + n + 1)\Omega$. Notice that the
characteristic eigenvalue parameter, i.e.~the eigenvalue responsible for the bifurcation, 
is the chemical potential. 
Thus, given an initial energy $\mu$ at the linear limit, we increase the chemical potential (and therefore the number of atoms $N$) 
in order to approach to the strongly nonlinear regime that can lead to the discovery of new states originating from secondary bifurcations.

For the numerical characterization of the stability (and possible classification of the instability) for each state, we refer to 
\cite{Charalampidis_et_al2018}.

\subsection{Finite Element discretization}\label{sec:FOM_GP}

In this subsection, we apply the method presented in Sec.~\ref{sec:FOM} to
problem \eqref{eq:N_Schrodinger}. We recall that the solution $\phi$ to eq.~\eqref{eq:N_Schrodinger}
is a complex function. Let $\varphi$ and $\psi$ be its real and imaginary part, respectively. 
Let us introduce:
\begin{gather*}
a(X, Y) = \frac{1}{2}\int_D{\nabla X \cdot \nabla Y} d\bx , \quad
b(X, Y) = \frac{1}{2}\Omega\int_D{|\boldsymbol{r}|^2 X \cdot Y}d\bx , \\
d(X, Y; \mu) = \mu\int_D{X \cdot Y}d\bx , \quad 
n(Y, Z) = \int_D{|Z|^2 Z \cdot Y}d\bx , \\
c(X, Y, Z) = \int_D{[2 (X \cdot Z)Z + |Z|^2 X] \cdot Y}d\bx , 
\end{gather*}

The generic $k$-th iteration of the Newton's method \eqref{linearnewtgal} reads: 
seek $\delta X_h = (\delta \varphi_h, \delta \psi_h) \in V_h$  with $V_h \subset (H_0^1)^2$, such that 
\begin{align}
a(\delta X_h, Y_h) + b(\delta X_h, Y_h) - d(\delta X_h, Y_h; \mu) + c(\delta X_h, X_h^k, Y_h) = \cl
a(X_h^k, Y_h) + b(X_h^k, Y_h) - d(X_h^k, Y_h; \mu) + n(Y_h, X_h) \quad   \forall\,  Y_h \in  V_h. \label{eq:weakformgalerkin}
\end{align}
Note that for simplicity of notation, we have not specified that the solution $X_h$ depends on 
parameter $\mu$. 

Let us introduce the finite element discretization matrices:
\begin{equation}
\begin{aligned}
(\mathbb{A}_h)_{ij} &= a(E^j, E^i) \ ,  \quad
(\mathbb{B}_h)_{ij} = b(E^j, E^i) \ , \\
(\mathbb{C}_h)_{ij} &= c(E^j, X_h^k, E^i) \ , \quad
(\mathbb{D}_h)_{ij} = c(E^j, E^i) \ .
\label{eq:matricesdef}
\end{aligned}
\end{equation}
The $\R^{N_h \times N_h}$ Jacobian matrix $\mathbb{J}(\vec{X}_h^k(\mu); \mu)$ can
be written as
\begin{equation*}
\mathbb{J}(\vec{X}_h^k(\mu); \mu) = \mathbb{A}_h + \mathbb{B}_h - \mathbb{D}_h + \mathbb{C}_h.
\end{equation*}

Next, we will apply the ROM technique presented in Sec.~\ref{sec:ROM}.

\subsection{Application of the RB method}

We present the problem that has to be solved in the online phase 
of the RB method applied to the Gross--Pitaevskii equation.
Given $\mu \in \D$, at every iteration $k$ of Newton's method
\eqref{eq:weakrednewton} we seek $\delta X_N = (\delta \phi_N, \delta \psi_N) \in V_N$ such that 
\begin{align}
a(\delta X_N, Y_N) + & b(\delta X_N, Y_N) - d(\delta X_N, Y_N; \mu) + c(\delta X_N, Y_N, X_N^K) = \cl 
& a(X_N^k, Y_N) + b(X_N^k, Y_N) - d(X_N^k, Y_N; \mu) + n(Y_N, X_N^k) \label{eq:weakformreduced}
\end{align}
$\forall\,  Y_N \in  V_N$. 

The reduced Jacobian $\mathbb{J}_N(\vec{X}_N^k(\mu); \mu) \in \R^{N\times N}$ can be written as
\begin{equation}
\label{matrixformred}
\mathbb{J}_N(\vec{X}_N^k(\mu); \mu)  = \mathbb{A}_N + \mathbb{B}_N - \mu\mathbb{D}_N + \mathbb{C}_N
\end{equation}
where:
\begin{gather*}
\mathbb{A}_N = \mathbb{V}^T\mathbb{A}_h\mathbb{V} \ ,  \quad
\mathbb{B}_N = \mathbb{V}^T\mathbb{B}_h\mathbb{V} \ ,  \\
\mathbb{D}_N = \mathbb{V}^T\mathbb{D}_h\mathbb{V} \ ,  \quad
\mathbb{C}_N = \sum_{n=1}^{N}X_N^{(n)}\mathbb{V}^T\mathbb{C}_h(\Sigma^n)\mathbb{V}
\label{eq:matricesdefred} 
\end{gather*} 
are the reduced matrices written in terms of finite element matrices \eqref{eq:matricesdef}
and transformation matrix \eqref{eq:matrix_V}. 

In general, the time saving promised by the online-offline strategy are enabled by 
the so called \textit{affine decomposition} \cite{hesthaven2015certified}, which makes the computations in the online phase independent form the (usually very high) number of degrees of freedom $N_h$. 
Here, however, reduced matrix $\mathbb{C}_N$ introduced above depends on $\mu$ through the solution computed 
at each step of Newton's method. 
Thus, we will need an affine-recovery technique called Empirical Interpolation Method (EIM) \cite{barrault04:_empir_inter_method}
in order to obtain substantial savings of the computational time during the online phase, 
as demonstrated in Sec.~\ref{sec:hyper}.

\section{Results}\label{sec:num_res}

The high fidelity (or full order) approximations used for the results presented in this section
were computed with FEniCS \cite{fenics,LoggMardalEtAl2012a,AlnaesBlechta2015a}, 
while we used RBniCS \cite{rbnics} for the reduced order approximations.

\subsection{Validation of the full order method}\label{sec:validation}

To validate the Full Order Method for the Gross--Pitaevskii
as described in Sec.~\ref{sec:FOM_GP}, we consider a test
proposed in  \cite{Charalampidis_et_al2018}. We approximate the solution to
eq.~\eqref{eq:N_Schrodinger} in domain $D =(-12,12)^2$, with 
homogeneous Dirichlet boundary conditions on the entire boundary of $D$.
We recall that we set $\Omega = 0.2$ and $\D = [0, 1.2]$. For the space discretization, we use 
$\mathbb{P}_2$ finite elements and a mesh with 6889 elements. 
 
\begin{figure}[h]
\begin{center}
\begin{overpic}[scale=.42,grid=false]{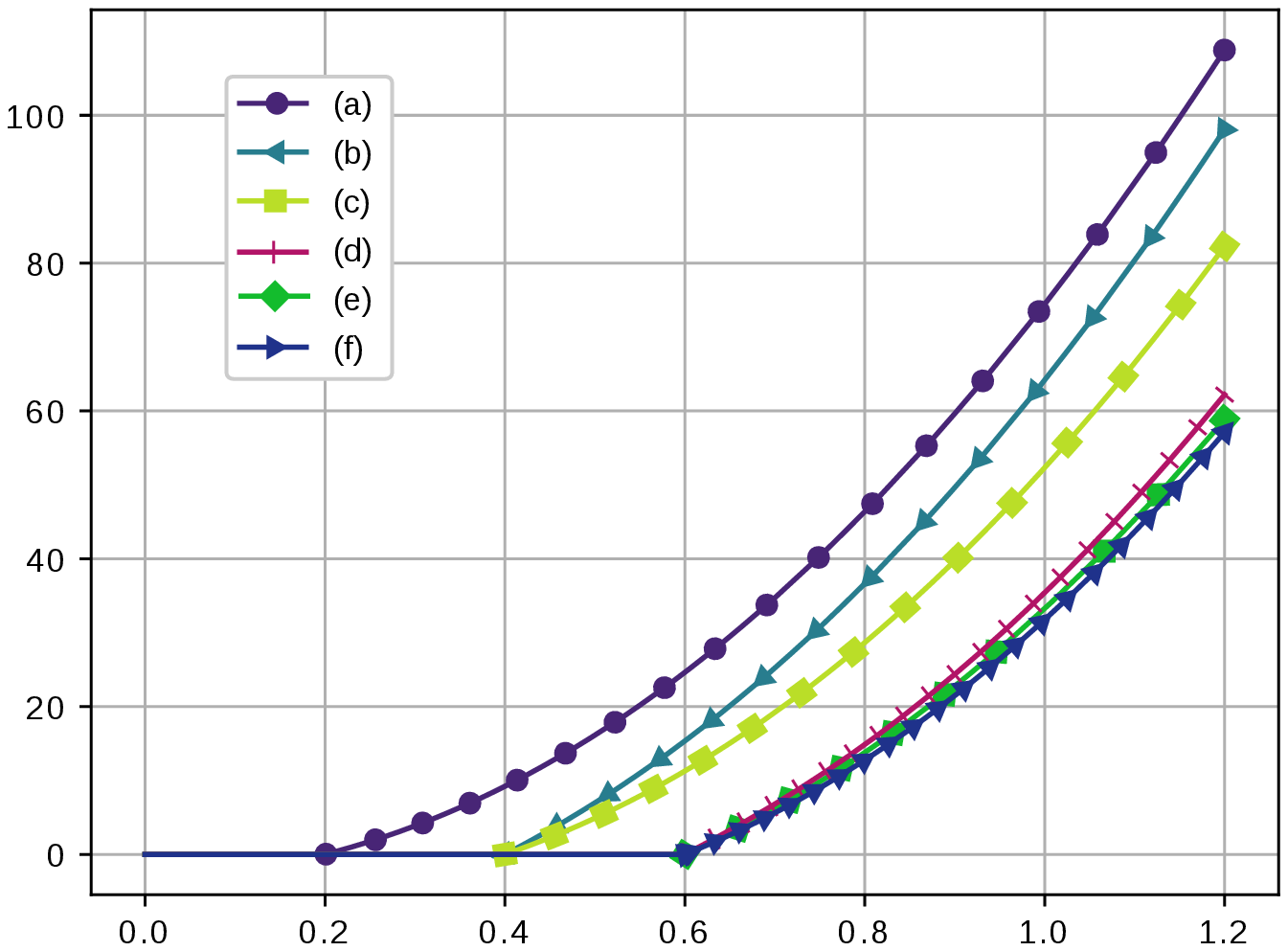}
    \put(52,0){\scriptsize $\mu$}
    \put(-3,43){\scriptsize $N$}
  \end{overpic}
  \hspace{2cm}
  \begin{overpic}[scale=.42,grid=false]{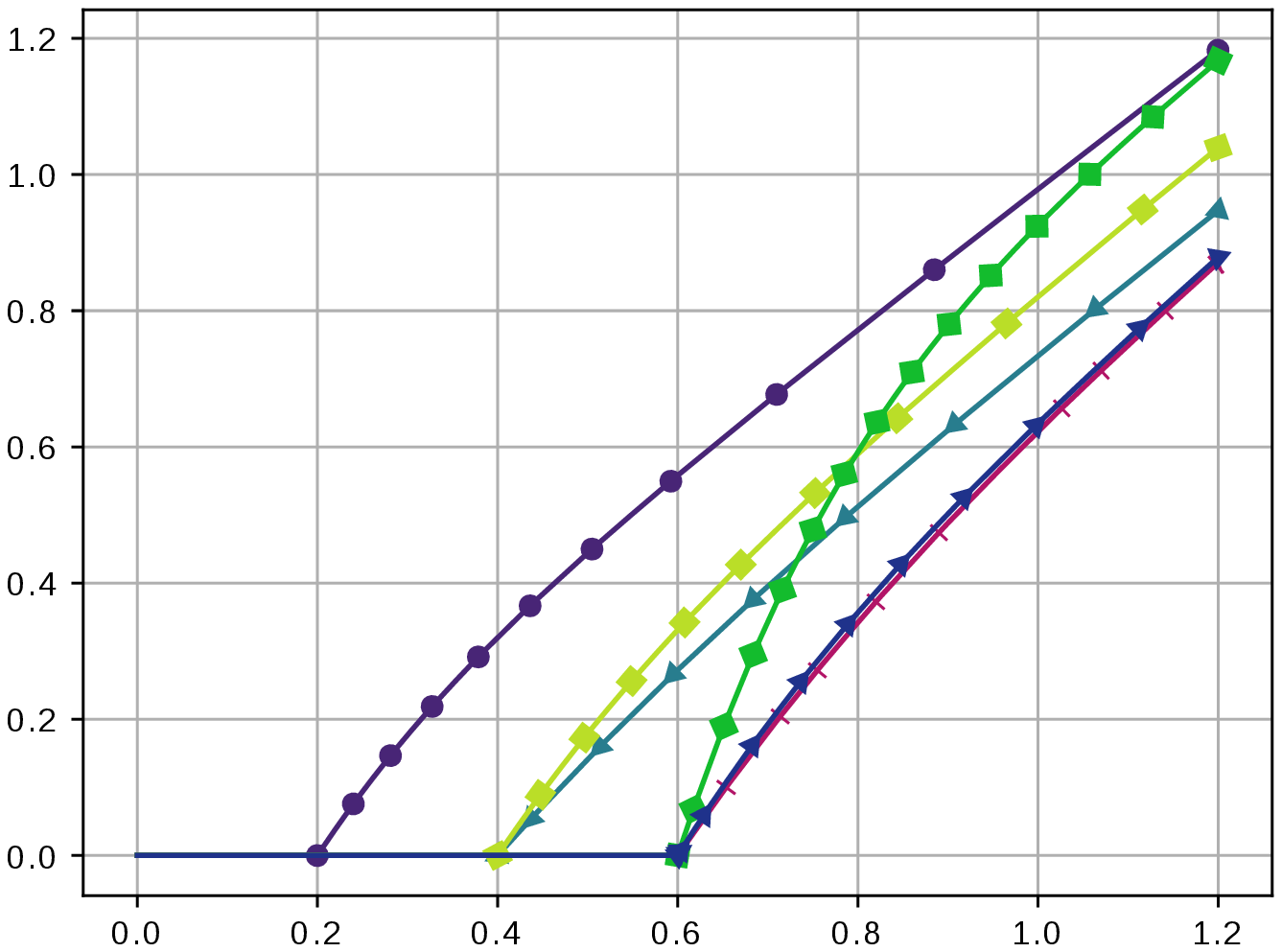}
    \put(52,0){\scriptsize $\mu$}
    \put(-7,45){\scriptsize $|| \rho ||_\infty$}
  \end{overpic}
\caption{Bifurcation diagram obtained with the Full Order Method: number of bosons $N$ (left) 
and infinity norm of density $\rho$ (right) plotted against the chemical potential $\mu$. The labels (a)-(f) are referred to solutions in Fig.~\ref{fig:rho_6branches}
}\label{fig:bif_FOM}
\end{center}
\end{figure}

 Fig.~\ref{fig:bif_FOM} displays the FOM bifurcation diagram in the $\mu$-$N$ plane (left)
and in the $\mu$-$|| \rho||_\infty$ plane (right).

The full reconstruction of the bifurcation diagram requires a proper initialization of our algorithm, 
in particular for Newton's method. For this purpose, we rely on the linear limit of the system as specified in~\eqref{herm_guess}.
Indeed, to approximate each branch in Fig.~\ref{fig:bif_FOM}, we assign an initial guess proportional to the product
of Hermite polynomials $H_m$ and $H_n$, where $m$ and $n$ chosen according to the value of the critical point $\mu_{crit} = (m + n + 1)\Omega$. In the general case (i.e., no analytic information is available), one can recover such an initial guess either 
from a linearized eigenvalue problem \cite{pichirozza} or thorough a deflation method \cite{Charalampidis_et_al2018, pintore2019efficient}.

\begin{figure}[h]
	\centering
	\begin{subfigure}{.31\textwidth}
		\centering
		\includegraphics[width=\textwidth]{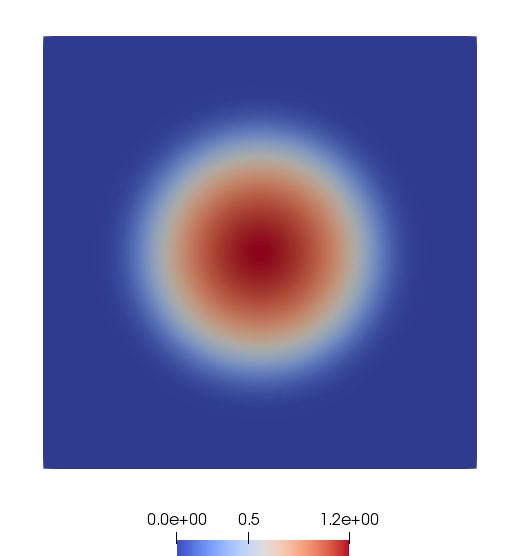}
		\caption{ground state}\label{fig:rho_branch1}
	\end{subfigure}
	\begin{subfigure}{.31\textwidth}
		\centering
		\includegraphics[width=\textwidth]{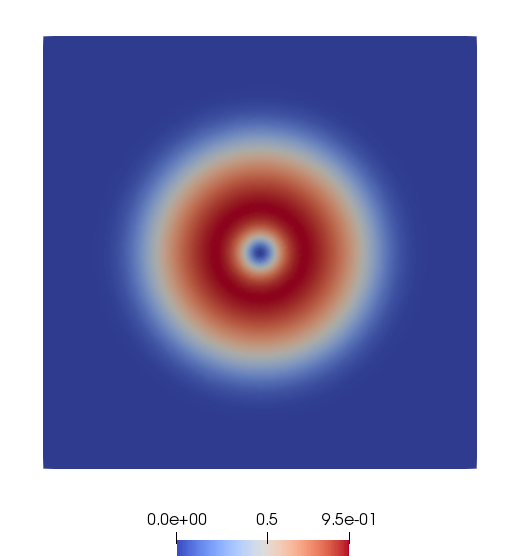}
		\caption{single charge vortex}\label{fig:rho_branch2a}
	\end{subfigure}
	\begin{subfigure}{.31\textwidth}
		\centering
		\includegraphics[width=\textwidth]{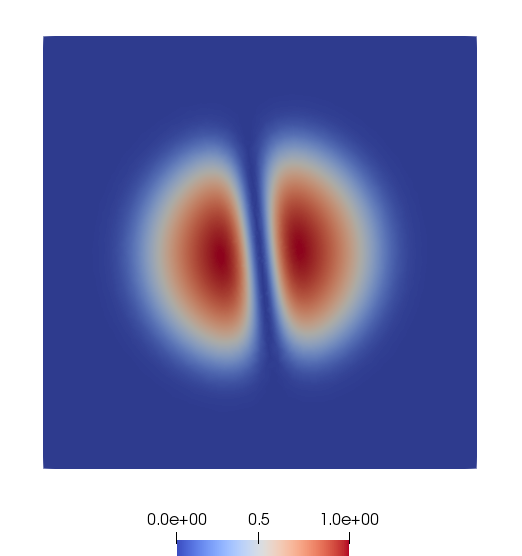}
		\caption{1-dark soliton stripe}\label{fig:rho_branch2b}
	\end{subfigure}
		\begin{subfigure}{.31\textwidth}
		\centering
		\includegraphics[width=\textwidth]{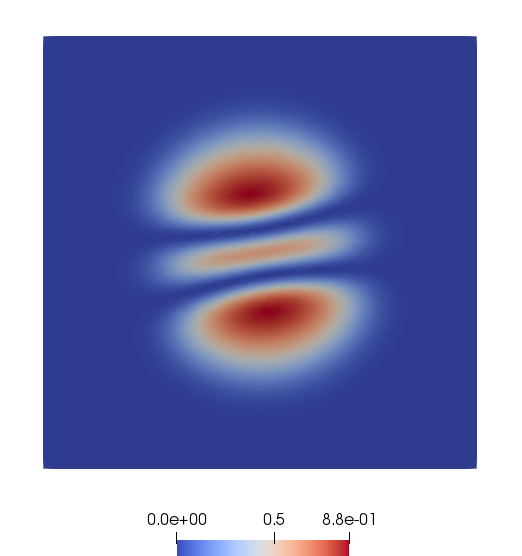}
		\caption{2-dark soliton stripe}\label{fig:rho_branch3c}
	\end{subfigure}
	\begin{subfigure}{.31\textwidth}
		\centering
		\includegraphics[width=\textwidth]{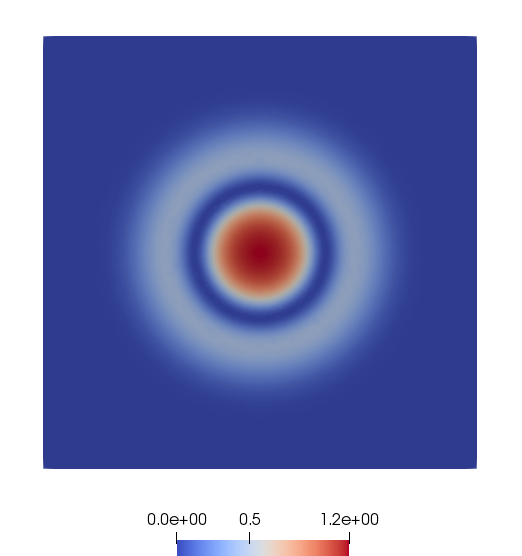}
		\caption{ring dark soliton}\label{fig:rho_branch3b}
	\end{subfigure}
	\begin{subfigure}{.31\textwidth}
		\centering
		\includegraphics[width=\textwidth]{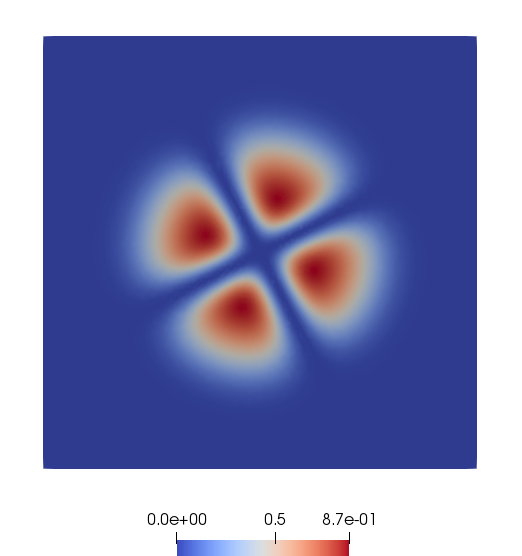}
		\caption{dark soliton cross}\label{fig:rho_branch3a}		
	\end{subfigure}
	\caption{Density functions computed with the full order method 
	for $\mu = 1.2$. Each plot is associated to one of the 6 solution branches in Fig.~\ref{fig:bif_FOM}.
	From (a) to (f) the number of bosons $N$ is decreasing.}
	\label{fig:rho_6branches}
\end{figure}

These diagrams show the first three bifurcation points and the relative
non-uniqueness of the solution with respect to the parameter $\mu$.  
As $\mu$ is increased, the sequence of events is as follows. 
The ground state $|0, 0\rangle$ is the system simplest state. Its linear eigenfunction $\phi_{0,0}$ has
corresponding eigenvalue $\mu = \Omega$. The ground state is generically stable, thus
no further bifurcations occur from this state \cite{doi:10.1137/1.9781611973945}. 
As expected, a unique solution branch departs from $\mu = \Omega$
in Fig.~\ref{fig:bif_FOM}. A representative density function for this branch is shown in Fig.~\ref{fig:rho_branch1}.
 We see no further bifurcation for $\Omega \leq \mu < 2\Omega$.
The first interesting events in terms of bifurcation analysis
occur for $\mu = 2 \Omega$, with $n+m = 1$: two
branches, associated to $|0, 1\rangle$ 
and $|1, 0\rangle$, bifurcate from point $(2 \Omega,0)$ in the $\mu$-$N$ and $\mu$-$|| \rho||_\infty$ planes 
\cite{Middelkamp_et_al2010,CONTRERAS2016265}. 
Indeed, from point $(2 \Omega,0)$ in Fig.~\ref{fig:bif_FOM} we observe the two expected branches.
Representative density functions for these two branches are reported in Fig.~\ref{fig:rho_branch2a} and \ref{fig:rho_branch2b}.
The next, more complicated, case of bifurcations emanates from  point $(3 \Omega,0)$, with $n+m = 2$.
In Fig.~\ref{fig:bif_FOM}, we see that three branches depart from this point, associated to 
$|1, 1\rangle$, $|0, 2\rangle$, and $|2, 0\rangle$. The corresponding representative
densities are shown in Fig.~\ref{fig:rho_branch3a}, \ref{fig:rho_branch3b}, and \ref{fig:rho_branch3c}.
Finally, all the points without marker in Fig.~\ref{fig:bif_FOM} correspond to the 
non-physical solution $\phi = 0$ that exists since there are no external forces in eq.~\eqref{eq:N_Schrodinger}.


Fig.~\ref{fig:rho_6branches} displays the  density functions 
associated to $\mu = 1.2$ and all the 6 solution branches in Fig.~\ref{fig:bif_FOM}. 
We observe the richness of density patterns in order of decreasing $N$. 
In particular, we see the ground state $|0, 0\rangle$ in Fig.~\ref{fig:rho_branch1}, 
the single charge vortex $|0, 1\rangle$ in Fig.~\ref{fig:rho_branch2a},
the 1-dark soliton stripe $|1, 0\rangle$ in Fig.~\ref{fig:rho_branch2b},
the dark soliton cross $|1, 1\rangle$ in Fig.~\ref{fig:rho_branch3a},
the ring dark soliton $|0, 2\rangle$  in Fig.~\ref{fig:rho_branch3b},
the 2-dark soliton stripe $|2, 0\rangle$ in Fig.~\ref{fig:rho_branch3c}.
Notice that the 6 branches Fig.~\ref{fig:bif_FOM} are related to the first three eigenvalues.
For example, the second bifurcation stems from a double eigenvalue and thus we have two branches. 
This phenomenon is called multiple bifurcations. 
The stability property of these branches are different for each case, i.e the single charge vortex is always stable while the 1-dark soliton stripe is subject to multiple secondary bifurcations. These properties can be easily studied using standard techniques 
(see, e.g., \cite{Charalampidis_et_al2018}). 
 The results in Fig.~\ref{fig:bif_FOM} and \ref{fig:rho_6branches} are in excellent
agreement with the results reported in  \cite{Charalampidis_et_al2018}, indicating that the
mesh that we use is sufficiently refined for this study.

The overall simulation time required to complete the diagrams
in Fig.~\ref{fig:bif_FOM} is roughly 96 minutes 
with continuation step $\Delta \mu = 1.25\cdot 10^{-3}$. 

\subsection{Tracing bifurcation diagrams with the reduced order approach}

\subsubsection{One-parameter study}

In this section, we present the results obtained with our Reduced Order Method as described in Sec.~\ref{sec:ROM}
and compare them with the FOM results reported in Sec.~\ref{sec:validation}. The only parameter that varies is $\mu$
in interval $[0, 1.2]$.

Concerning the construction of the reduced manifold, we employed a training set for the POD 
with cardinality $N_{train} = 160$ for each one of the six branches.
Setting the POD tolerance to $10^{-9}$, we obtain a global basis of dimension $N = 51$. 
In the online phase, we reconstruct the reduced bifurcation diagram for all the 961 equally spaced 
points in $\D = [0, 1.2]$ used in the high-fidelity bifurcation diagram shown in Fig.~\ref{fig:bif_FOM}.
Such points correspond to continuation step $\Delta \mu = 1.25\cdot 10^{-3}$.

\begin{figure}[h]
\begin{center}
\begin{overpic}[scale=.55,grid=false]{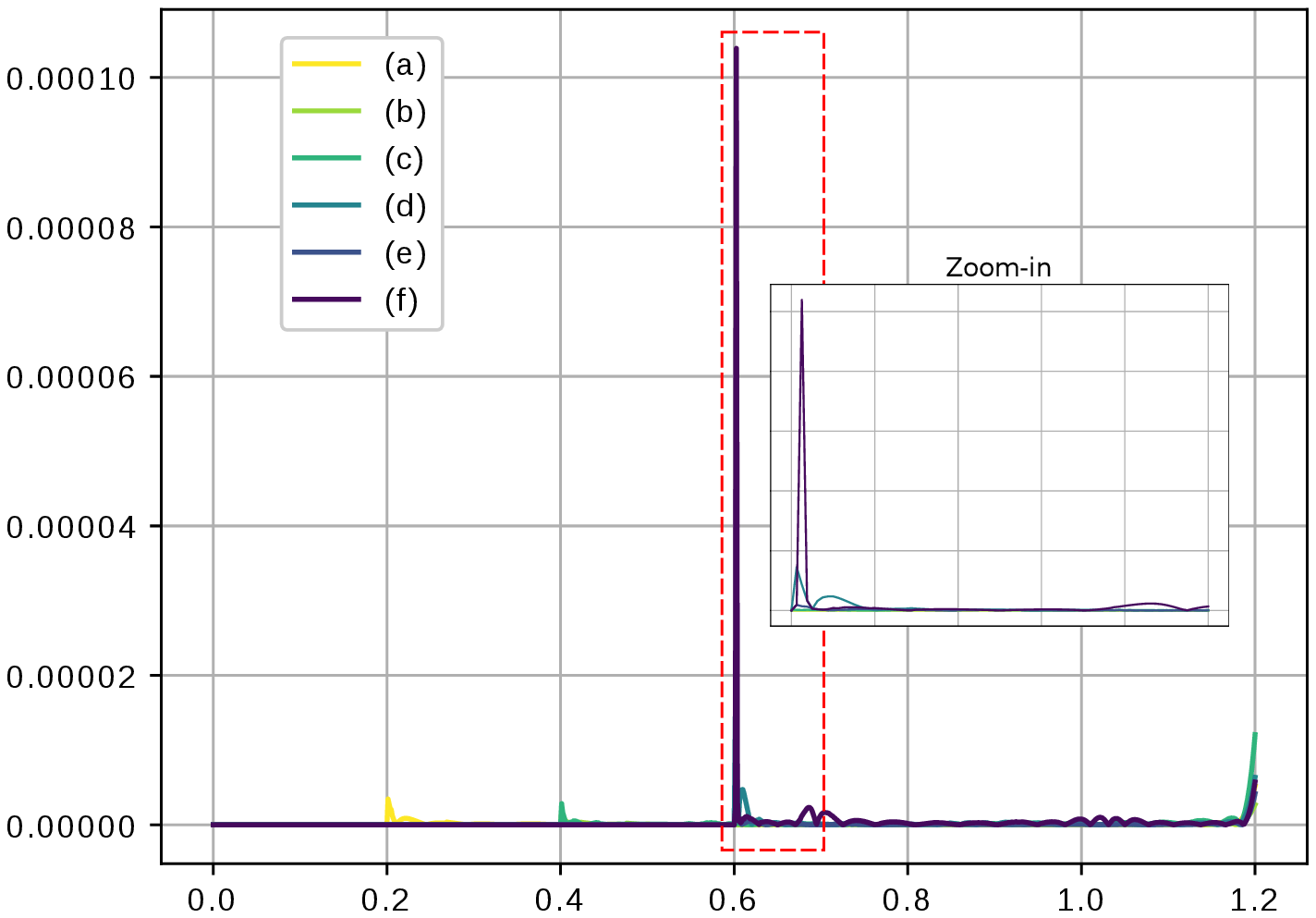}
    \put(55,2){\scriptsize $\mu$}
    \put(-6,40){\scriptsize $E_N$}
  \end{overpic}
  \hfill
  \begin{overpic}[scale=.55,grid=false]{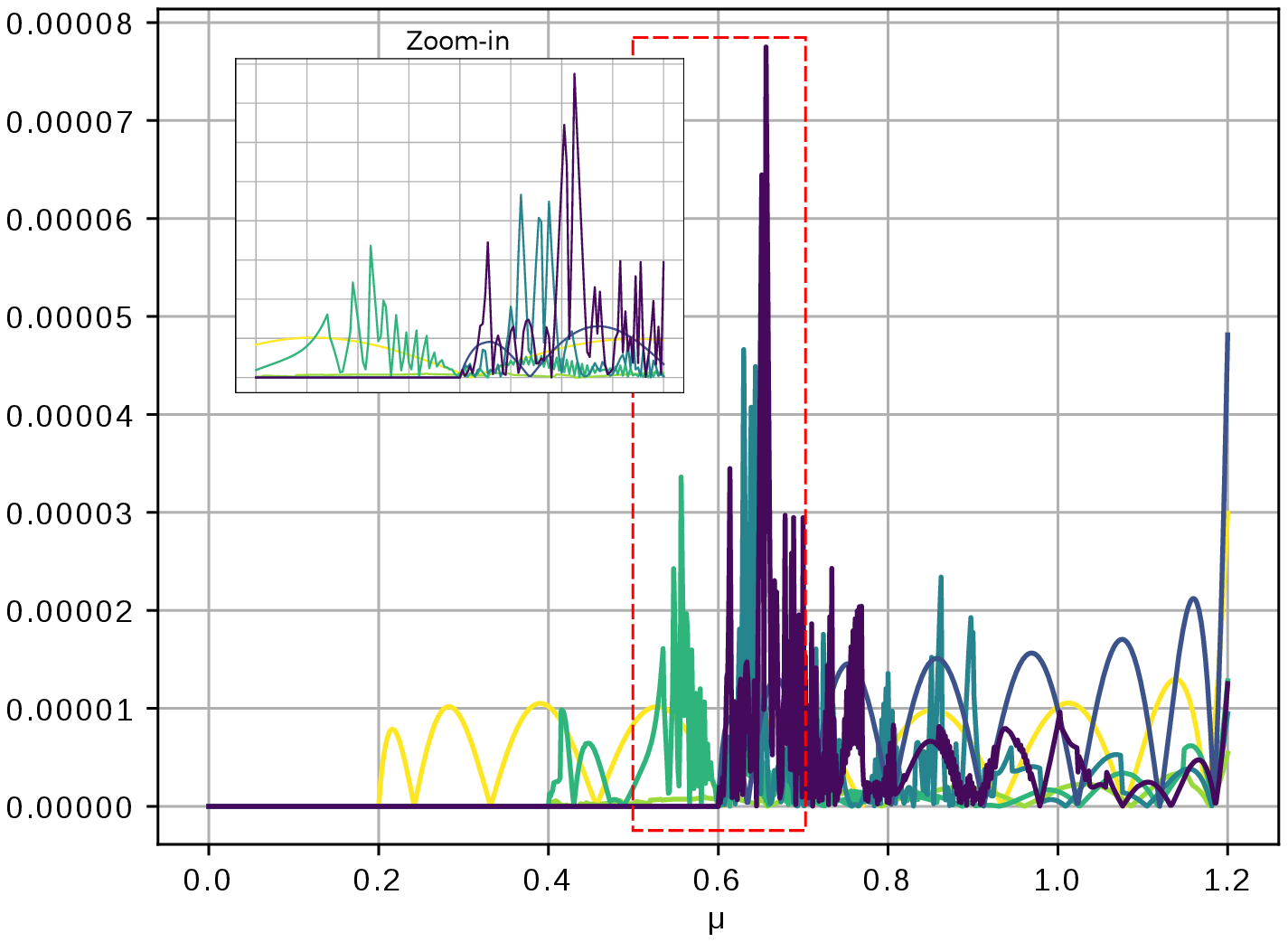}
   \put(-6,40){\scriptsize $E_\rho$}
  \end{overpic}
\caption{Difference in absolute value between the branches of the bifurcation diagram computed with FOM and ROM
in the $\mu$-$N$ plane (top) and in the $\mu$-$|| \rho||_\infty$ plane (bottom),
i.e.~reduced order errors $E_N$ (top) and $E_\rho$ (bottom) defined in \eqref{N_rho_error}. The labels (a)-(f) refer 
to the solutions reported in Fig.~\ref{fig:rho_6branches}.
}\label{fig:bif_ROM}
\end{center}
\end{figure}

Fig.~\ref{fig:bif_ROM} shows reduced order errors 
\begin{align}\label{N_rho_error}
E_N = | N_{h} - N_{N} | \quad \text{and} \quad E_\rho = | || \rho_{h}||_{\infty} - ||\rho_{N}||_{\infty}, 
\end{align}
i.e.~the difference in absolute value between the branches of the bifurcation diagram computed with FOM and ROM
in the $\mu$-$N$ plane (top) and in the $\mu$-$|| \rho||_\infty$ plane (bottom).
In Fig.~\ref{fig:bif_ROM}, we see that the largest peaks are associated to the $|1, 1\rangle$
branch at $\mu = 0.6$. 
In general, it is expected to have larger errors at the bifurcation points 
where differentiability with respect to the parameter $\mu$ is lost. We infer that 
the errors are largest at $\mu = 0.6$ due to the more complicated solution structure
(compare Fig.~\ref{fig:rho_branch3a} to the other panels in Fig.~\ref{fig:rho_6branches}).
From Fig.~\ref{fig:bif_ROM} (bottom), we see the largest error is of the order of $10^{-4}$ also for
the infinity norm of the density. However, we observe larger errors over interval $[0.2, 1.2]$ for $\mu$, 
as opposed to localized at $\mu = 0.6$ as in Fig.~\ref{fig:bif_ROM} (top).

\begin{figure}
	\centering
	\begin{subfigure}{.45\textwidth}
		\centering
		\includegraphics[width=\textwidth]{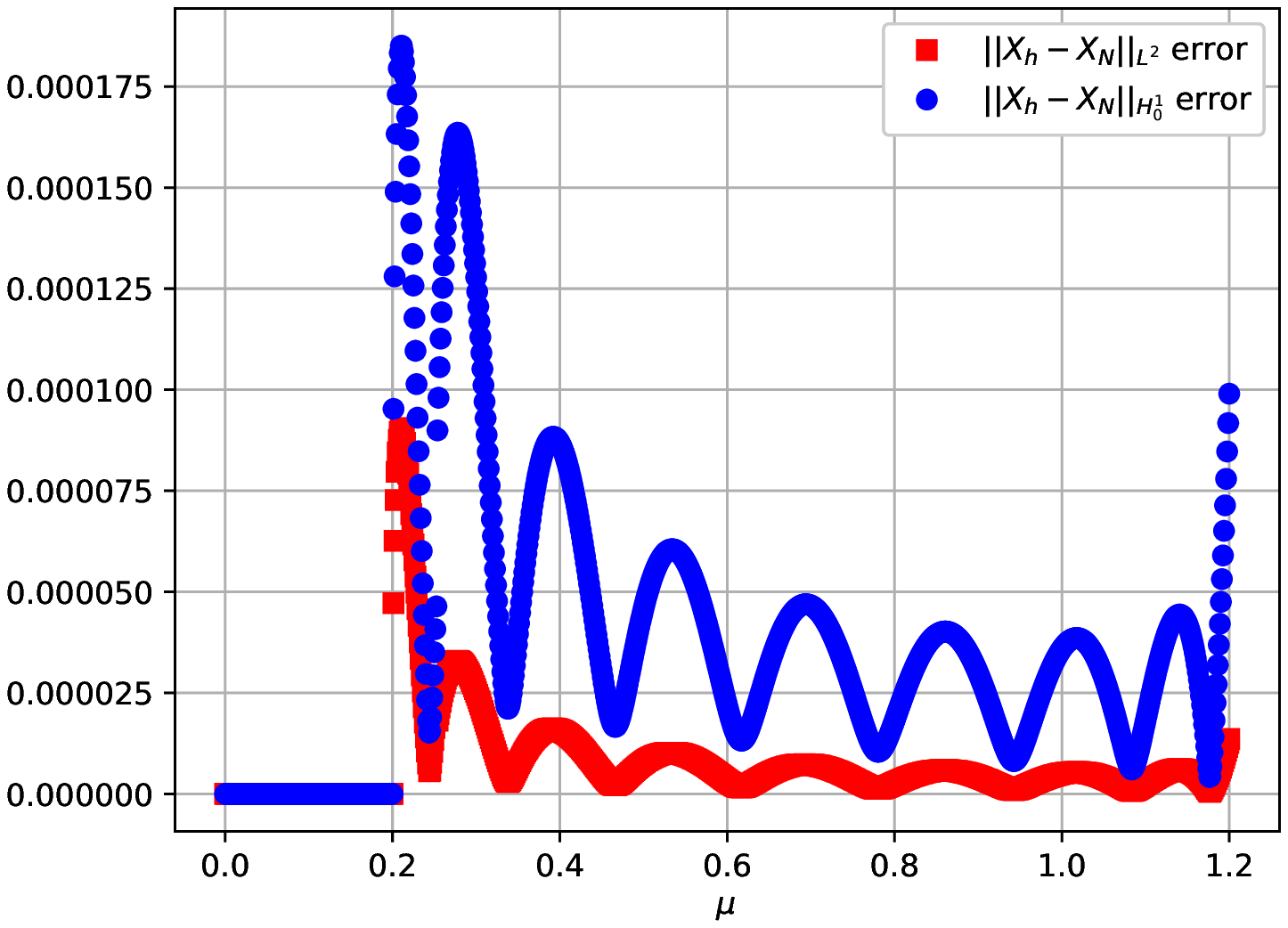}
		\caption{branch $|0,0\rangle$}\label{fig:rho_branch1_err}
	\end{subfigure}
	\begin{subfigure}{.45\textwidth}
		\centering
		\includegraphics[width=\textwidth]{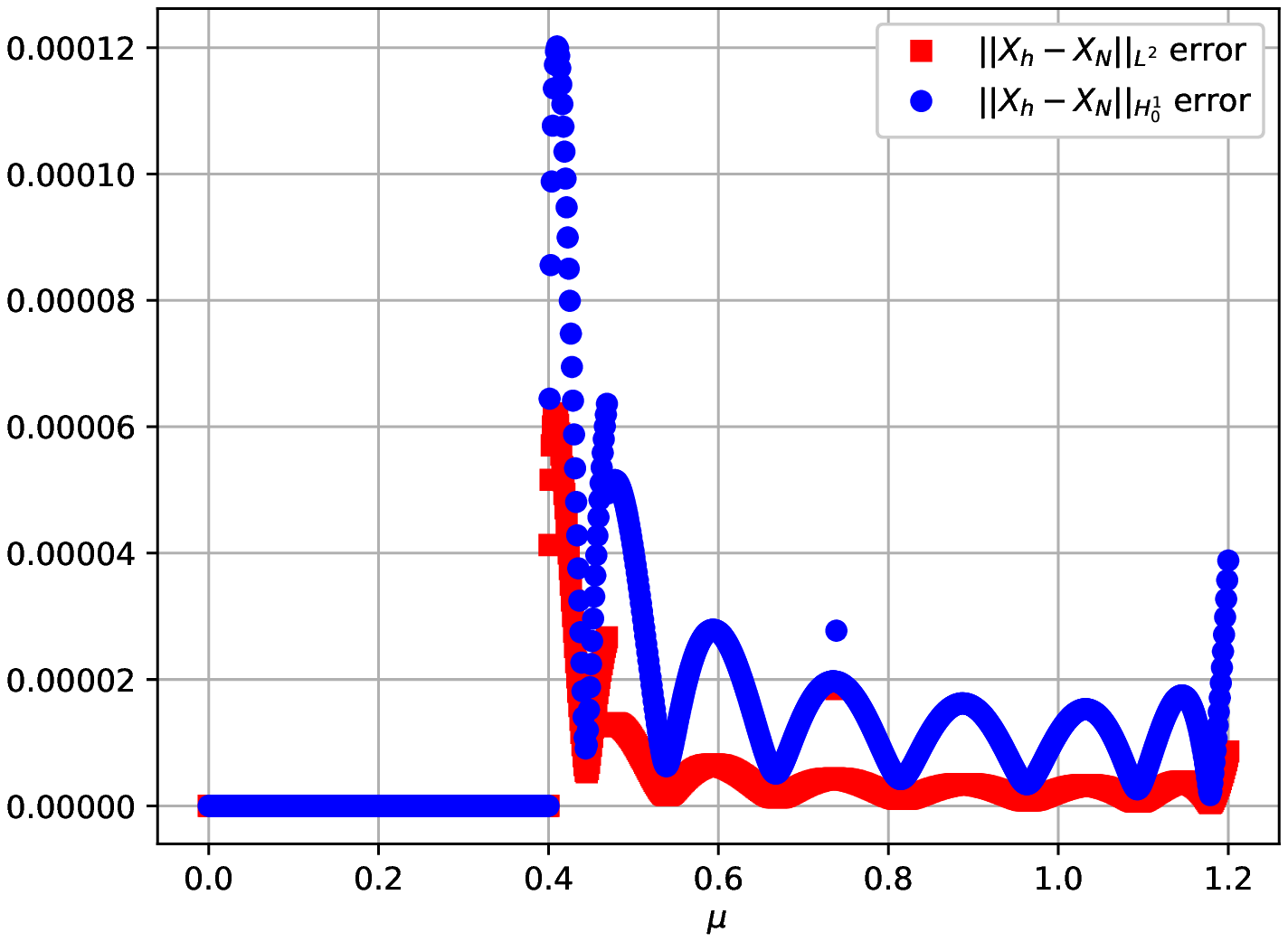}
		\caption{branch $|0,1\rangle$}\label{fig:rho_branch2a_err}
	\end{subfigure}
		\vspace{2em}
	
	\begin{subfigure}{.45\textwidth}
		\centering
		\includegraphics[width=\textwidth]{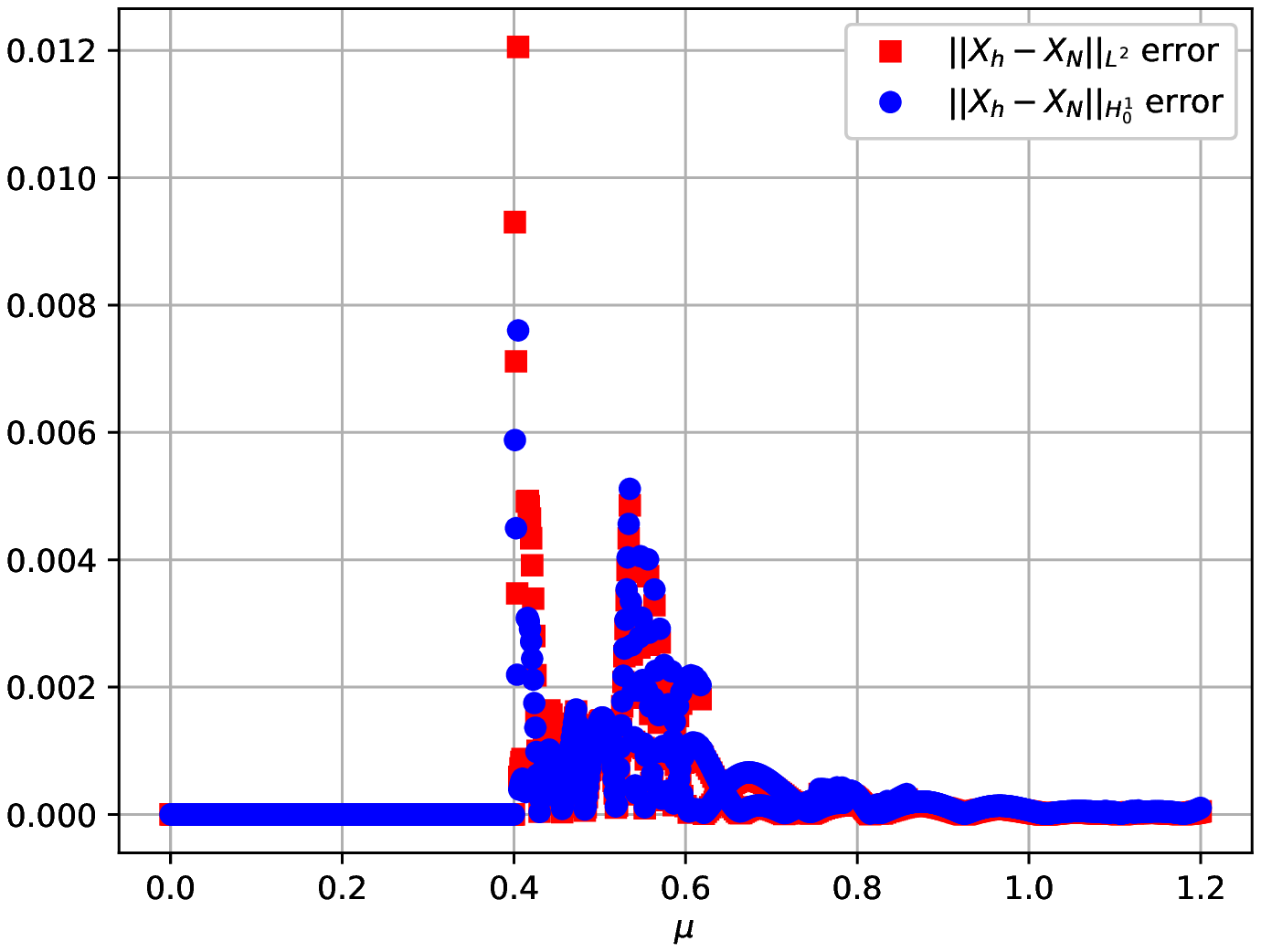}
		\caption{branch $|1,0\rangle$}\label{fig:rho_branch2b_err}
	\end{subfigure}
		\begin{subfigure}{.45\textwidth}
		\centering
		\includegraphics[width=\textwidth]{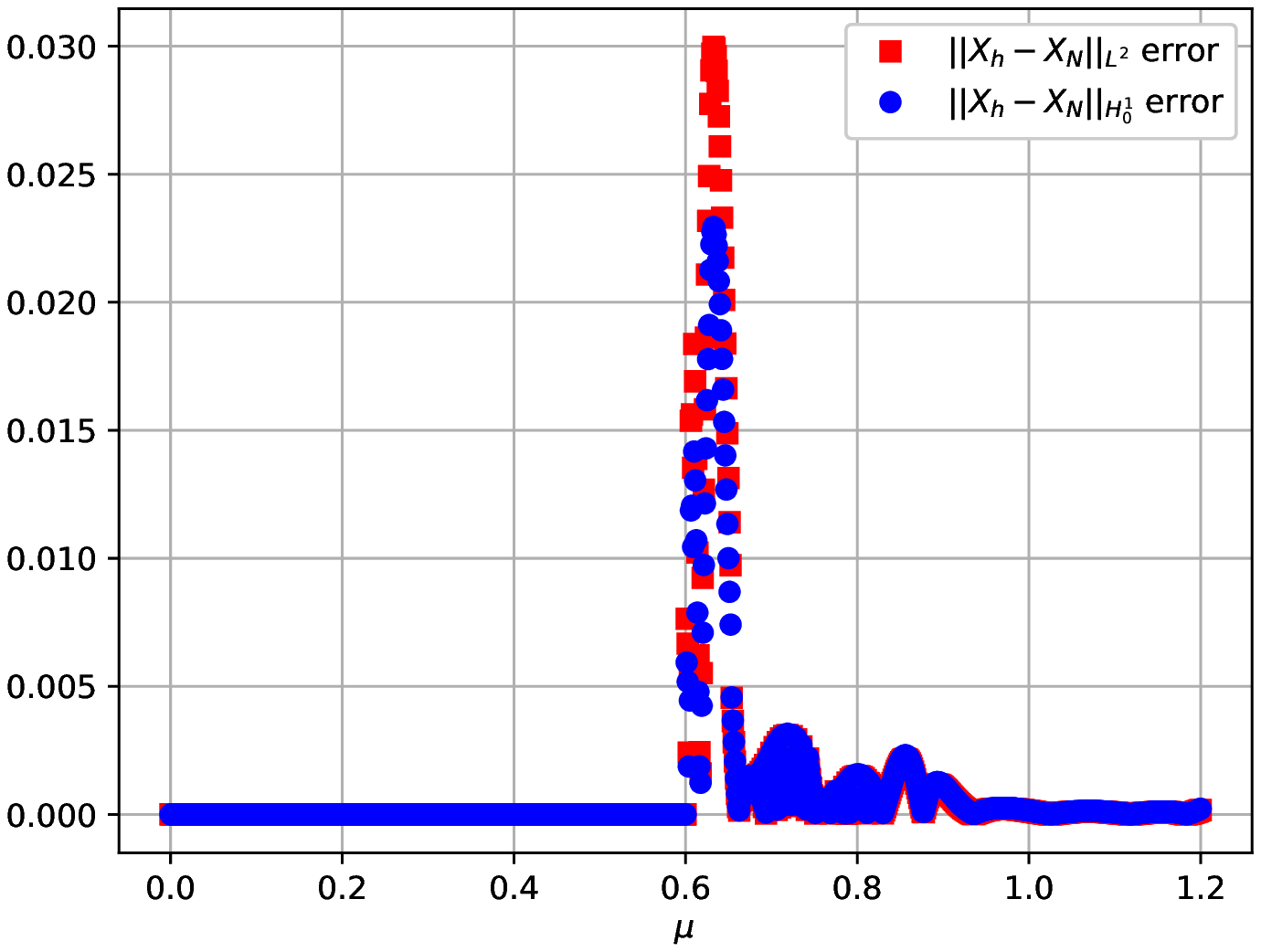}
		\caption{branch $|1,1\rangle$}\label{fig:rho_branch3a_err}
	\end{subfigure}
		\vspace{2em}
	
	\begin{subfigure}{.45\textwidth}
		\centering
		\includegraphics[width=\textwidth]{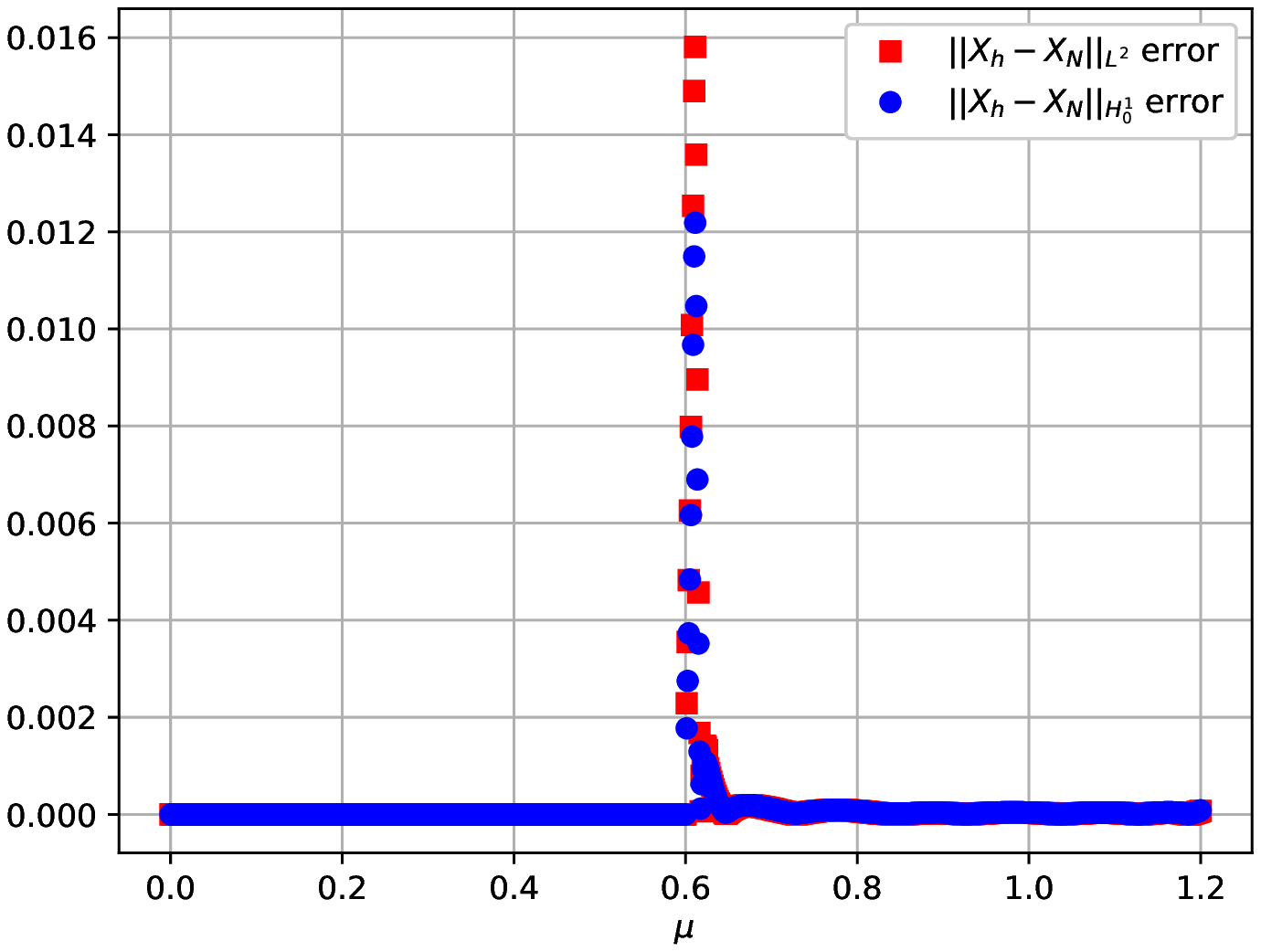}
		\caption{branch $|0,2\rangle$}\label{fig:rho_branch3b_err}
	\end{subfigure}
	\begin{subfigure}{.45\textwidth}
		\centering
		\includegraphics[width=\textwidth]{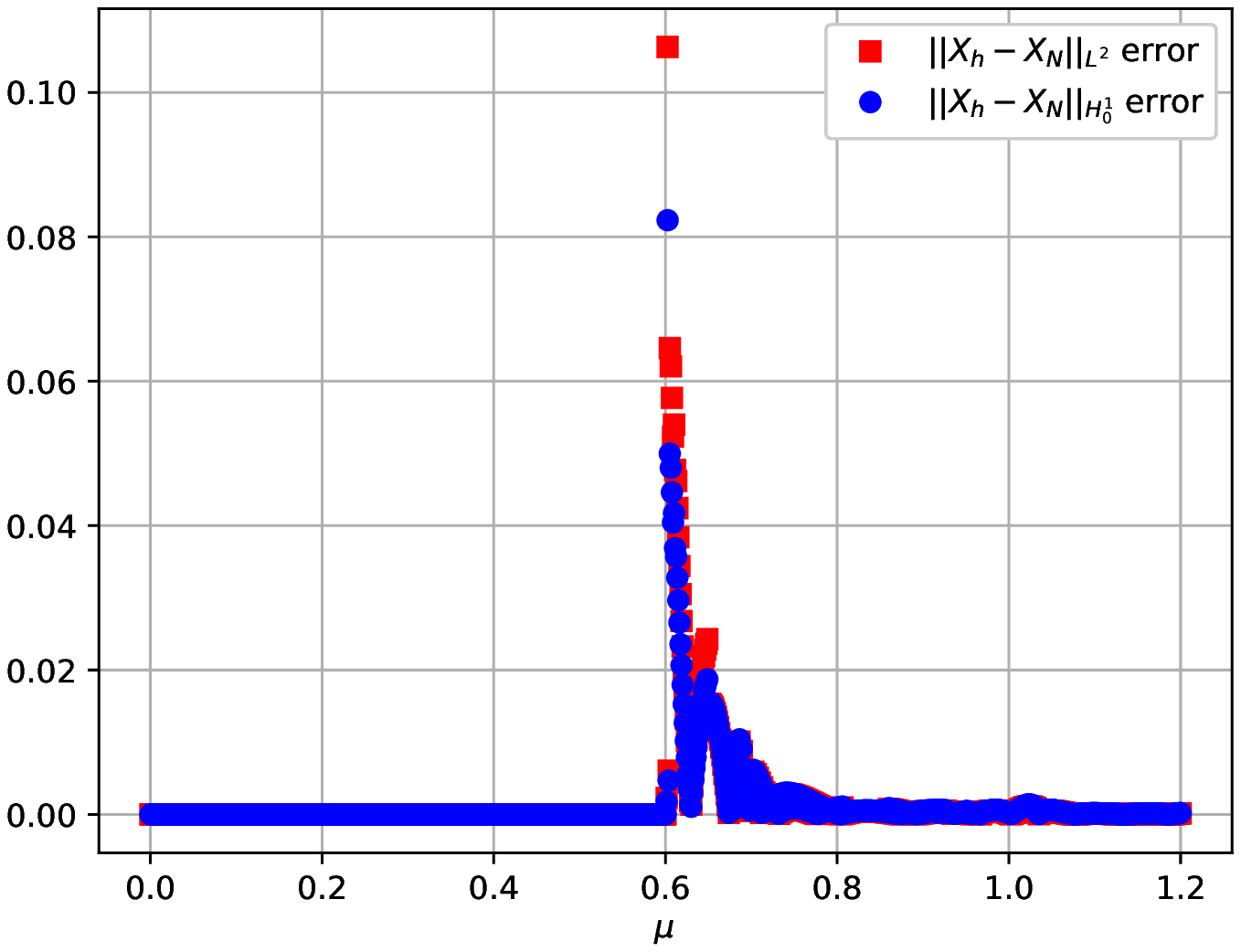}
		\caption{branch $|2,0\rangle$}\label{fig:rho_branch3c_err}
	\end{subfigure}
	\caption{Difference between $X_h$ (FOM solution) and $X_N$ (ROM solution) in the $L^2$ and $H^1_0$ norms for each of the six solution branches.}
\label{fig:RB_errors}
\end{figure}

As further evidence of the accuracy of our ROM approach, we plot in Fig.~\ref{fig:RB_errors}
the difference between $X_h$ (FOM solution) and $X_N$ (ROM solution) in the $L^2$ and $H^1_0$ norms. 
Recall that $X_h$ and $X_N$ consists of both real and imaginary part of the computed solution of
the Gross--Pitaevskii equation. 
As expected, for all the branches the largest errors in both norms
occur at the point where each branch departs from the horizontal axis, i.e. $\mu = 0.2$ for branch $|0, 0\rangle$, 
$\mu = 0.4$ for branches $|1, 0\rangle$ and $|0, 1\rangle$, and $\mu = 0.6$ for branches $|1, 1\rangle$, $|0, 2\rangle$,
and $|2, 0\rangle$. In addition, just like in Fig.~\ref{fig:bif_ROM} the largest errors are associated with branch $|1, 1\rangle$.

\begin{figure}
	\centering
	\begin{subfigure}{.45\textwidth}
		\centering
		\includegraphics[width=\textwidth]{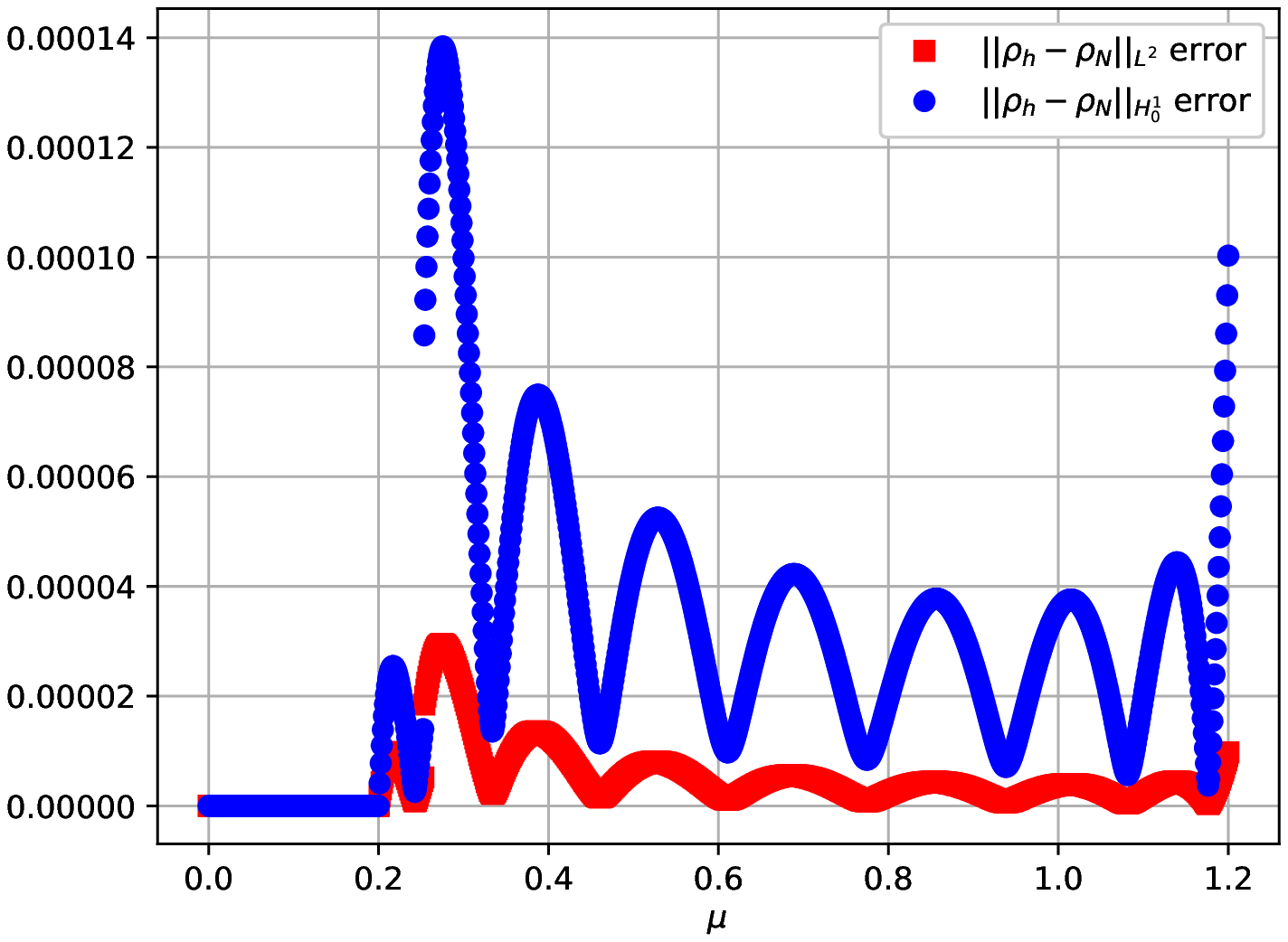}
		\caption{branch $|0,0\rangle$}\label{fig:rho_branch1_err_den}
	\end{subfigure}
	\begin{subfigure}{.45\textwidth}
		\centering
		\includegraphics[width=\textwidth]{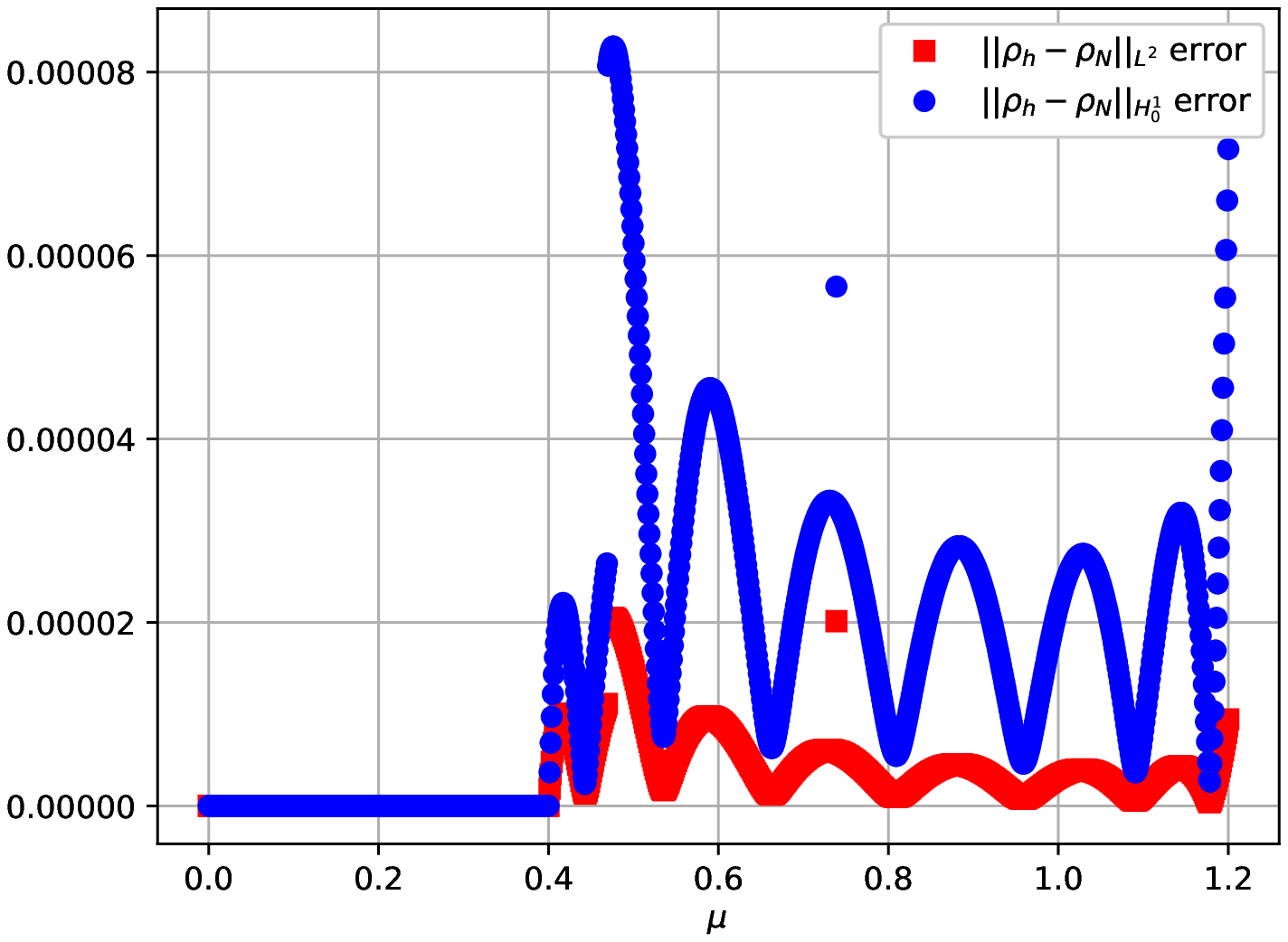}
		\caption{branch $|0,1\rangle$}\label{fig:rho_branch2a_err_den}
	\end{subfigure}
	\vspace{2em}
	
	\begin{subfigure}{.45\textwidth}
		\centering
		\includegraphics[width=\textwidth]{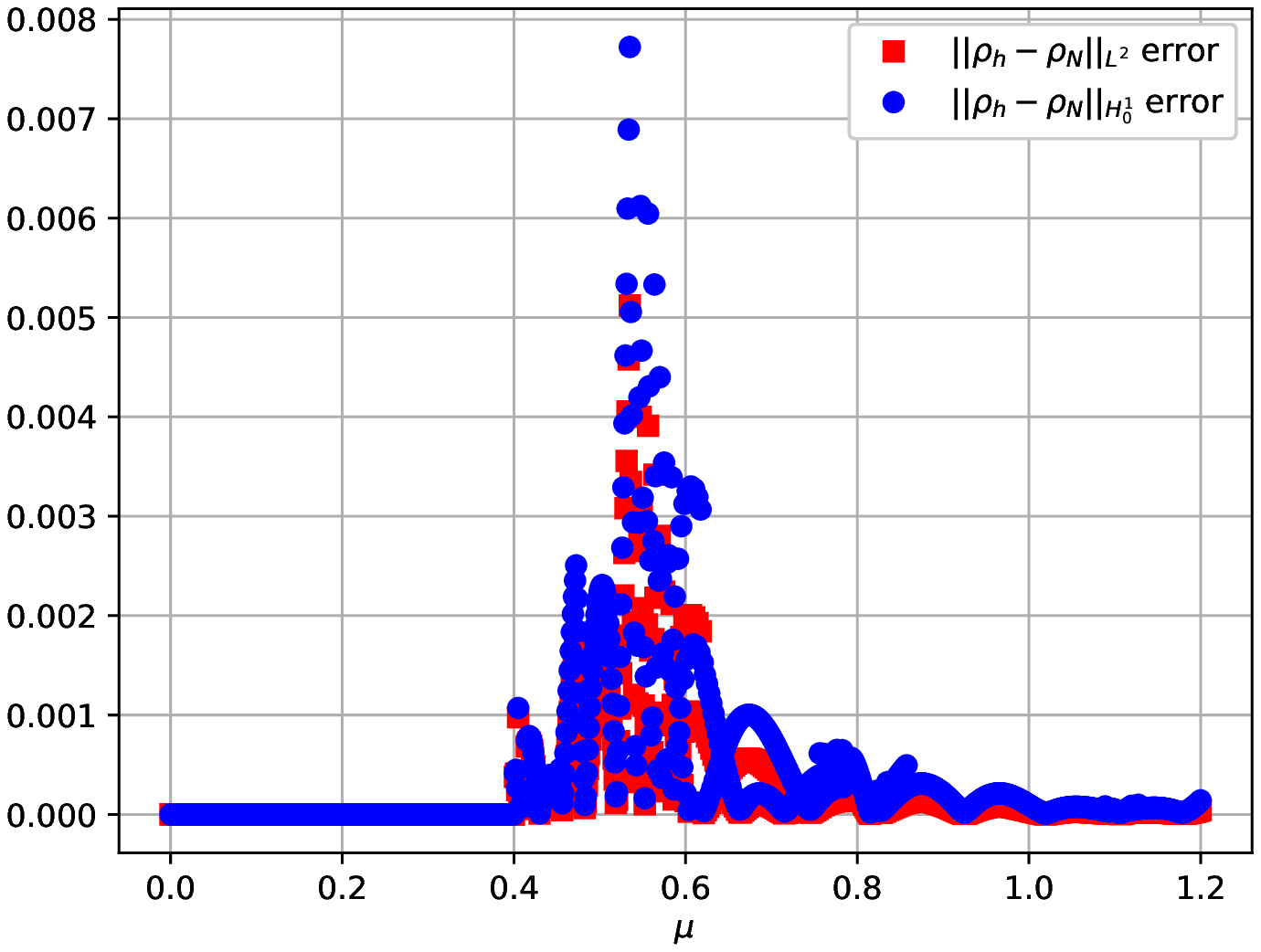}
		\caption{branch $|1,0\rangle$}\label{fig:rho_branch2b_err_den}
	\end{subfigure}
		\begin{subfigure}{.45\textwidth}
		\centering
		\includegraphics[width=\textwidth]{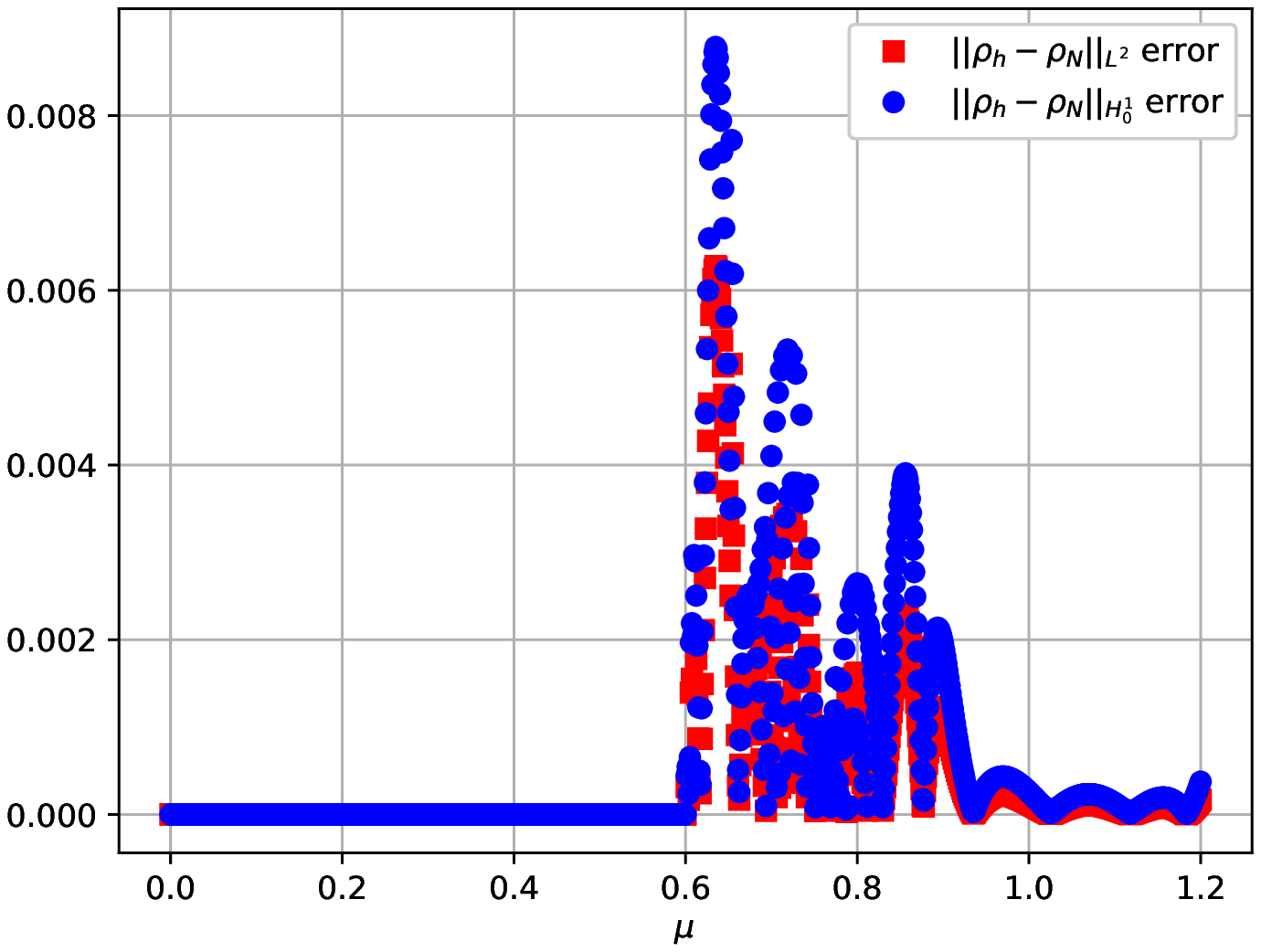}
		\caption{branch $|1,1\rangle$}\label{fig:rho_branch3a_err_den}
	\end{subfigure}
	\vspace{2em}
	
	\begin{subfigure}{.45\textwidth}
		\centering
		\includegraphics[width=\textwidth]{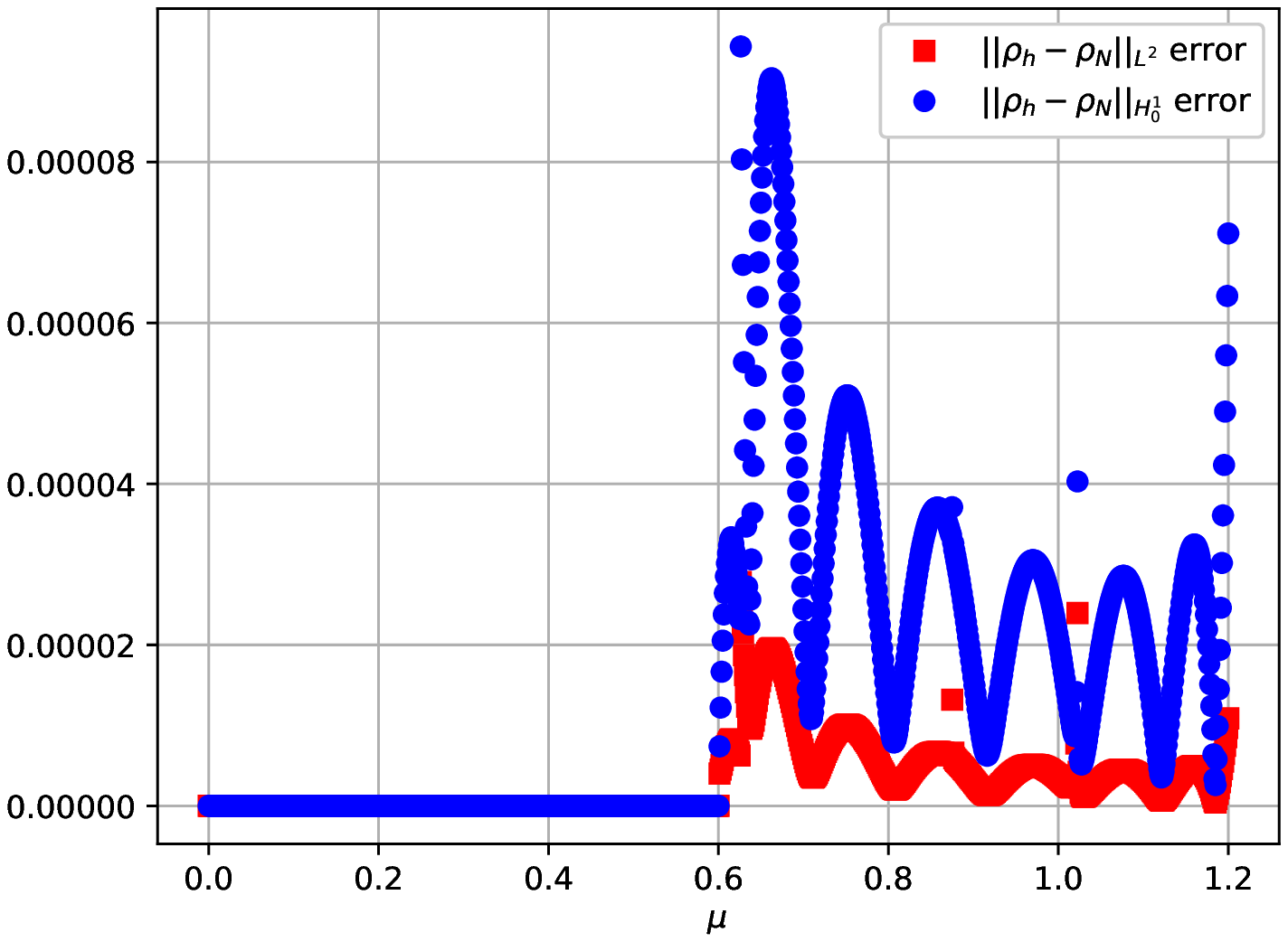}
		\caption{branch $|0,2\rangle$}\label{fig:rho_branch3b_err_den}
	\end{subfigure}
	\begin{subfigure}{.45\textwidth}
		\centering
		\includegraphics[width=\textwidth]{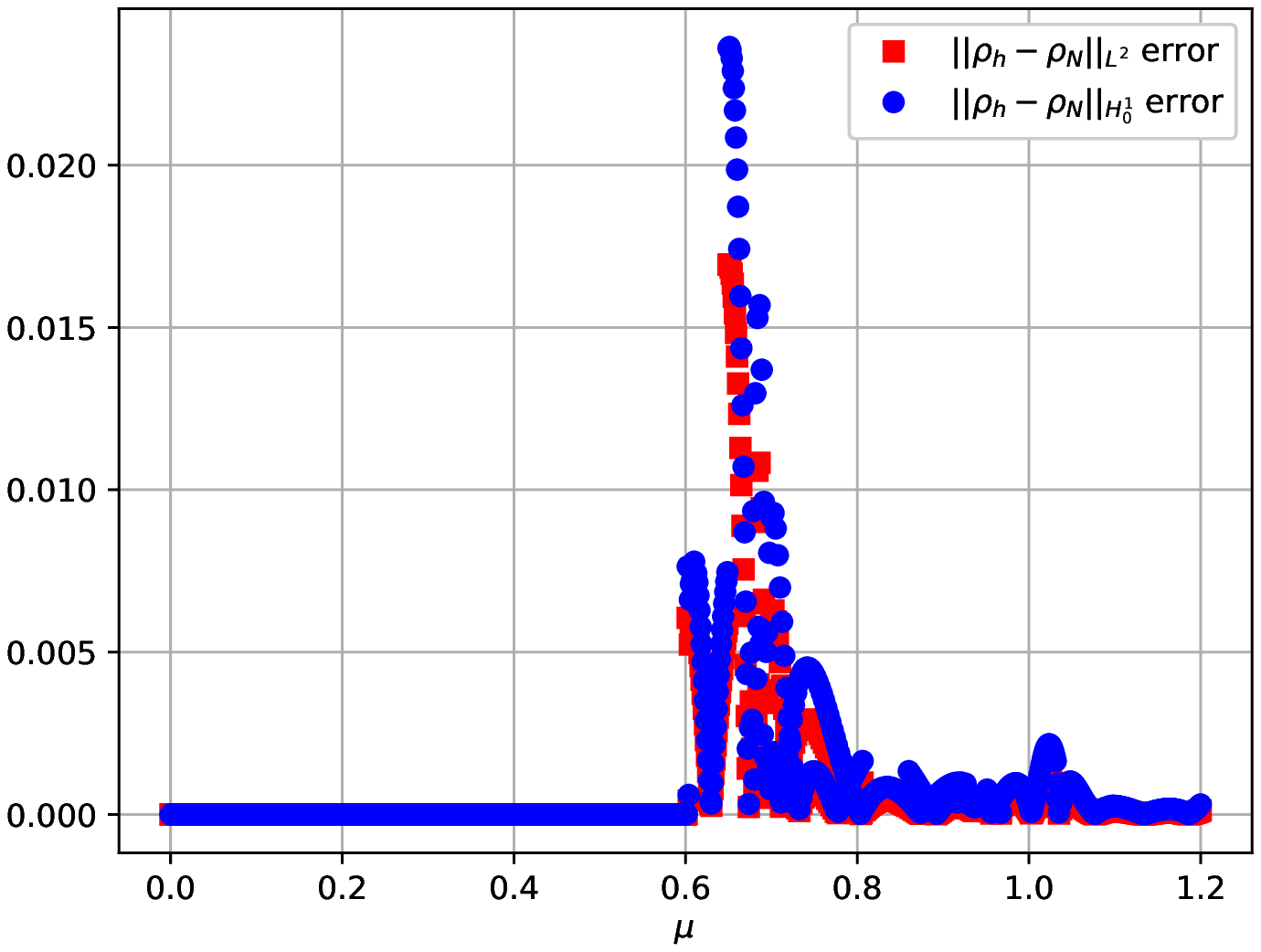}
		\caption{branch $|2,0\rangle$}\label{fig:rho_branch3c_err_den}
	\end{subfigure}
	\caption{Difference between $\rho_h$ (density computed with FOM) and $\rho_N$ (density computed with ROM) 
	in the $L^2$ and $H^1_0$ norms for each of the six solution branches.}
\label{fig:RB_errors_DEN}
\end{figure}

\begin{figure}
	\centering
	\begin{subfigure}{.32\textwidth}
		\centering
		\includegraphics[width=\textwidth]{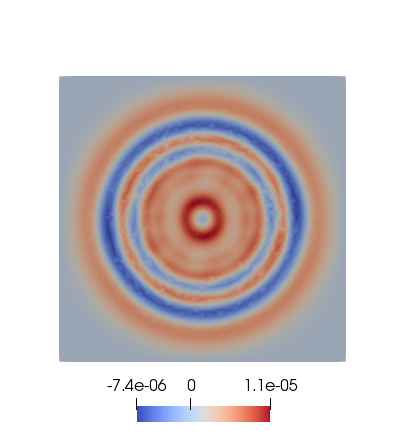}
		\caption{branch $|0,0\rangle$}\label{fig:rho_branch1_err_dens}
	\end{subfigure}
	\begin{subfigure}{.32\textwidth}
		\centering
		\includegraphics[width=\textwidth]{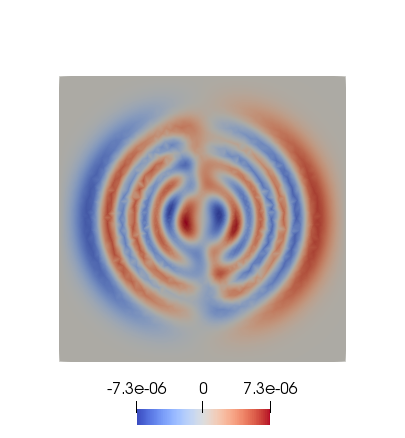}
		\caption{branch $|0,1\rangle$}\label{fig:rho_branch2a_err_dens}
	\end{subfigure}
	\begin{subfigure}{.32\textwidth}
		\centering
		\includegraphics[width=\textwidth]{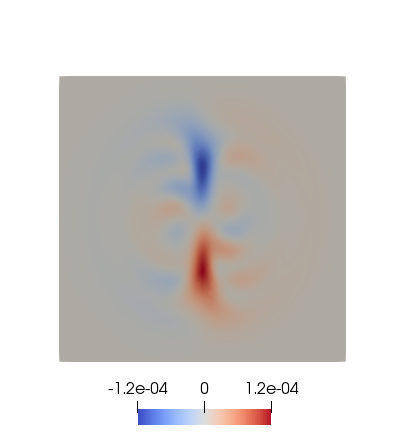}
		\caption{branch $|1,0\rangle$}\label{fig:rho_branch2b_err_dens}
	\end{subfigure}
		\begin{subfigure}{.32\textwidth}
		\centering
		\includegraphics[width=\textwidth]{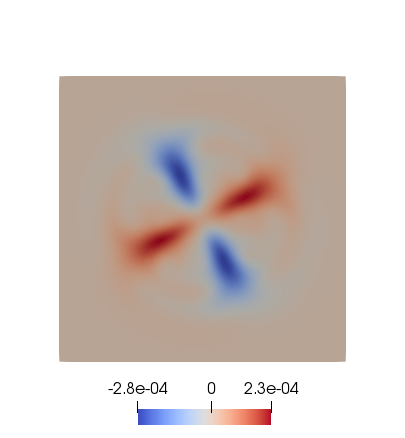}
		\caption{branch $|1,1\rangle$}\label{fig:rho_branch3a_err_dens}
	\end{subfigure}
	\begin{subfigure}{.32\textwidth}
		\centering
		\includegraphics[width=\textwidth]{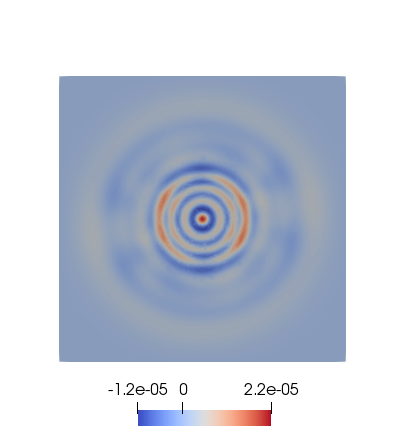}
		\caption{branch $|0,2\rangle$}\label{fig:rho_branch3b_err_dens}
	\end{subfigure}
	\begin{subfigure}{.32\textwidth}
		\centering
		\includegraphics[width=\textwidth]{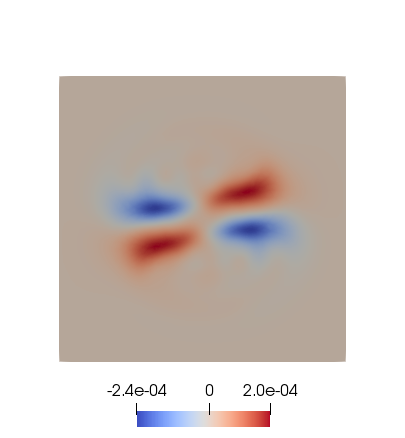}
		\caption{branch $|2,0\rangle$}\label{fig:rho_branch3c_err_dens}
	\end{subfigure}
	\caption{Difference between the density computed with FOM and ROM for the first six branches for $\mu = 1.2$.}
\label{fig:RB_errors_DENs}
\end{figure}

Fig.~\ref{fig:RB_errors_DEN} reports the difference between the density function $\rho$ computed with FOM and ROM 
in the $L^2$ and $H^1_0$ norms. In the case of the density, the largest errors for each branch occur for $\mu$
larger than the critical value where the branch departs from the horizontal axis. 
For a better understanding of how the ROM density compares with the FOM density,  
Fig.~\ref{fig:RB_errors_DENs} displays the difference $\rho_h - \rho_N$
for $\mu = 1.2$. 
We observe larger errors for the 1- and 2-dark soliton stripe and the dark soliton cross, i.e.~for the solutions
that do not have central symmetry. 

Fig.~\ref{fig:bif_ROM}-\ref{fig:RB_errors_DENs} show the ability of our ROM approach to accurately
reconstruct bifurcation diagrams as parameter $\mu$ varies. Clearly, it makes sense to set up the 
ROM machinery if there is a substantial gain in terms of computational time. 
Because of the nonlinearity in the Gross--Pitaevskii equation, which makes the computations
in the online phase dependent on the number of FEM degrees of freedom, the computational speed-up 
enabled by our ROM approach is only 1.1: it took 86 minutes 
to generate the data needed for the bifurcation diagrams, using 
continuation step $\Delta \mu = 1.25\cdot 10^{-3}$, accounting
only for the online phase computations.
Recall that it takes 96 minutes with FOM. These computational savings are not satisfactory,
especially if we were to include the cost for the offline phase. 
Before introducing an affine recovery technique to drastically improve computational efficiency,
we present a two-parameter study.

\subsubsection{Two-parameter study}\label{sec:two_parameter}

In this section, we plot the bifurcation diagram as
chemical potential $\mu$ varies in interval $[0, 1.2]$ and trap strength $\Omega$ varies in interval $[0.1, 0.3]$.

\begin{figure}
\begin{center}
\begin{overpic}[scale=.55,grid=false]{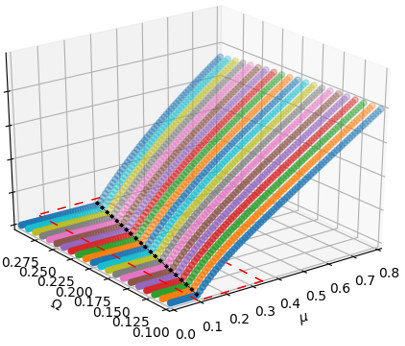}
    \put(97,45){\scriptsize $|| \rho ||_\infty$}
  \end{overpic}
\caption{First bifurcation in a two-parameter bifurcation diagram obtained with the Reduced Order Model: 
infinity norm of the density as chemical potential $\mu$ and trap strength $\Omega$
vary. The black dotted line in the highlighted red rectangle shows the critical value of $\mu$ for the first bifurcation: $\mu_{crit} = \Omega$.
}\label{fig:bif_3d_FOM_2}
\end{center}
\end{figure}

First, we focus on the first bifurcation, i.e.~branch $|0,0\rangle$.
Fig.~\ref{fig:bif_3d_FOM_2} show such branch in a two-parameter bifurcation diagram obtained with 
the Reduced Order Model as $\mu$ and $\Omega$ are varied. As expected from the theory, 
we see that as $\Omega$ increases the critical value of $\mu$ for the first bifurcation increases linearly.
Recall that $\mu_{crit} = \Omega$, for $m = n = 0$, and Fig.~\ref{fig:bif_3d_FOM_2} clearly shows it
(see the black dotted line in the highlighted red rectangle).
For Fig.~\ref{fig:bif_3d_FOM_2}, we used increment $\Delta \Omega = 0.01$ and 
continuation step $\Delta \mu = 1.25\cdot 10^{-3}$. The online phase computations
for the graph in Fig.~\ref{fig:bif_3d_FOM_2} took roughly 164 minutes. 
The corresponding time required by the FOM is approximately 249 minutes. 
So with our ROM approach we obtain a speed-up of 1.5, which represents an improvement
over the 1.1 speed-up for the one parameter study but it is still not enough to justify
the computational costs of the offline phase.

\begin{figure}
\begin{center}
\begin{overpic}[scale=.7,grid=false]{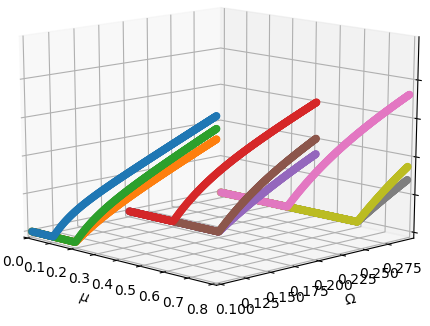}
    \put(100,45){\scriptsize $|| \rho ||_\infty$}
  \end{overpic}
\caption{First two bifurcations in a two-parameter bifurcation diagram obtained with the Reduced Order Model: 
infinity norm of the density as chemical potential $\mu$ varies for trap strength $\Omega = 0.1, 0.2, 0.3$. 
}\label{fig:bif_3d_FOM_1}
\end{center}
\end{figure}

Next, we focus on the first two bifurcations but restrict the attention to $\Omega = 0.1, 0.2, 0.3$. Fig.~\ref{fig:bif_3d_FOM_1}
shows the first three branches, i.e.~branches $|0,0\rangle$, $|0,1\rangle$, and
$|1,0\rangle$.
We see that our ROM approach successfully captures also the critical $\mu$ for the second 
bifurcation: branches $|0,1\rangle$ and $|1,0\rangle$ depart from $\mu_{crit} = 2\Omega$.
Just like for Fig.~\ref{fig:bif_3d_FOM_2}, we used continuation step $\Delta \mu = 1.25\cdot 10^{-3}$.
The time required by the online phase computations to plot the graph in  Fig.~\ref{fig:bif_3d_FOM_1} is
roughly 245 minutes. 

\subsection{Hyper-reduction techniques}\label{sec:hyper}

In this section, we use two affine-recovery techniques to reconstruct branch $|0,1\rangle$
for the one parameter study. These techniques are
called Empirical Interpolation Method (EIM) \cite{barrault04:_empir_inter_method} and Discrete Empirical Interpolation Method (DEIM) \cite{Chaturantabut2010} and they makes the 
online phase computations independent from the number of degrees of freedom of the finite element method.

As we have remarked earlier, if the affine dependence assumption is not fulfilled,
 the speedup of the online reconstruction can be limited, thereby compromising the whole methodology. 
With the hyper-reduction techniques one aims at approximating a general parametrized function 
$g: D \times \D \to \R$ by a sum of affine terms:
\begin{equation}\label{eq:hr1}
 g(\bx, \mu) \approx \texttt{I}[g_\mu](\bx) = \sum_{q=1}^Q{c_q(\mu)h_q(\bx)}
\end{equation}
where the basis functions $h_q$ are obtained by means of a linear combination of $Q$ snapshots $\{g_{\mu_q}\}_{q=1}^Q$ and the sample points are chosen through a Greedy approach. 
In \eqref{eq:hr1}, the coefficients of such expansion $c_q(\mu)$ are found by solving 
\begin{equation*}
\texttt{I}[g_\mu](\bx_j) = g_\mu(\bx_j)
\end{equation*}
in some particular points $\{\bx_j\}_{j=1}^Q$ of the domain $D$ called magic points. Hence, 
Empirical Interpolation strategies provide a discrete version of $\texttt{I}[g_\mu](x)$ as 
\begin{equation*}
g_M(\mu) = \mathbb{H}c(\mu), \qquad \qquad \mathbb{H} = \{h_q(\bx_j)\}_{(j,q)} \in \R^{Q \times Q}.
\end{equation*}

The main difference between EIM and DEIM is in the construction of $\mathbb{H}$. In fact, EIM embeds 
the construction of the basis inside the Greedy procedure, while DEIM exploits a POD on a set of snapshots. 
Moreover, the DEIM strategy starts with discretizing the nonlinearity, while EIM constructs 
the set of the magic points and the basis functions before the discretization step.

Let us express residual \eqref{linearnewtgal} as 
\begin{equation*}
G_h(\vec{X}_h; \mu) = \sum_{q=1}^{Q_g} \theta_q^g(\mu)G_h^q(\vec{X}_h),
\end{equation*}
where parameter affine dependency is guaranteed by the forms $\theta_q^g(\mu)$. The reduced residual \eqref{linearnewtred}
can then be expressed as:
\begin{equation*}
G_N(\vec{X}_N; \mu) = \sum_{q=1}^{Q_g} \theta_q^g(\mu)\mathbb{V}^TG_h^q(\mathbb{V}\vec{X}_N).
\end{equation*}
The hyper-reduction techniques described above provide the following affine approximation of the reduced residual
\begin{equation*}
G_N(\vec{X}_N; \mu) \approx \sum_{q=1}^{Q_g} \theta_q^g(\mu)c_q(\vec{X}_N; \mu)\mathbb{V}^T h^q
\end{equation*}
where $\{h^q\}_{q=1}^{Q_g}$ represent a suitable basis and $c_q$ the interpolation coefficients. The Jacobian matrix $\mathbb{J}_N(\vec{X}_N; \mu) $ is assembled in a similar fashion.

\begin{figure}[h]
\begin{center}
\begin{overpic}[width=.48\textwidth,grid=false]{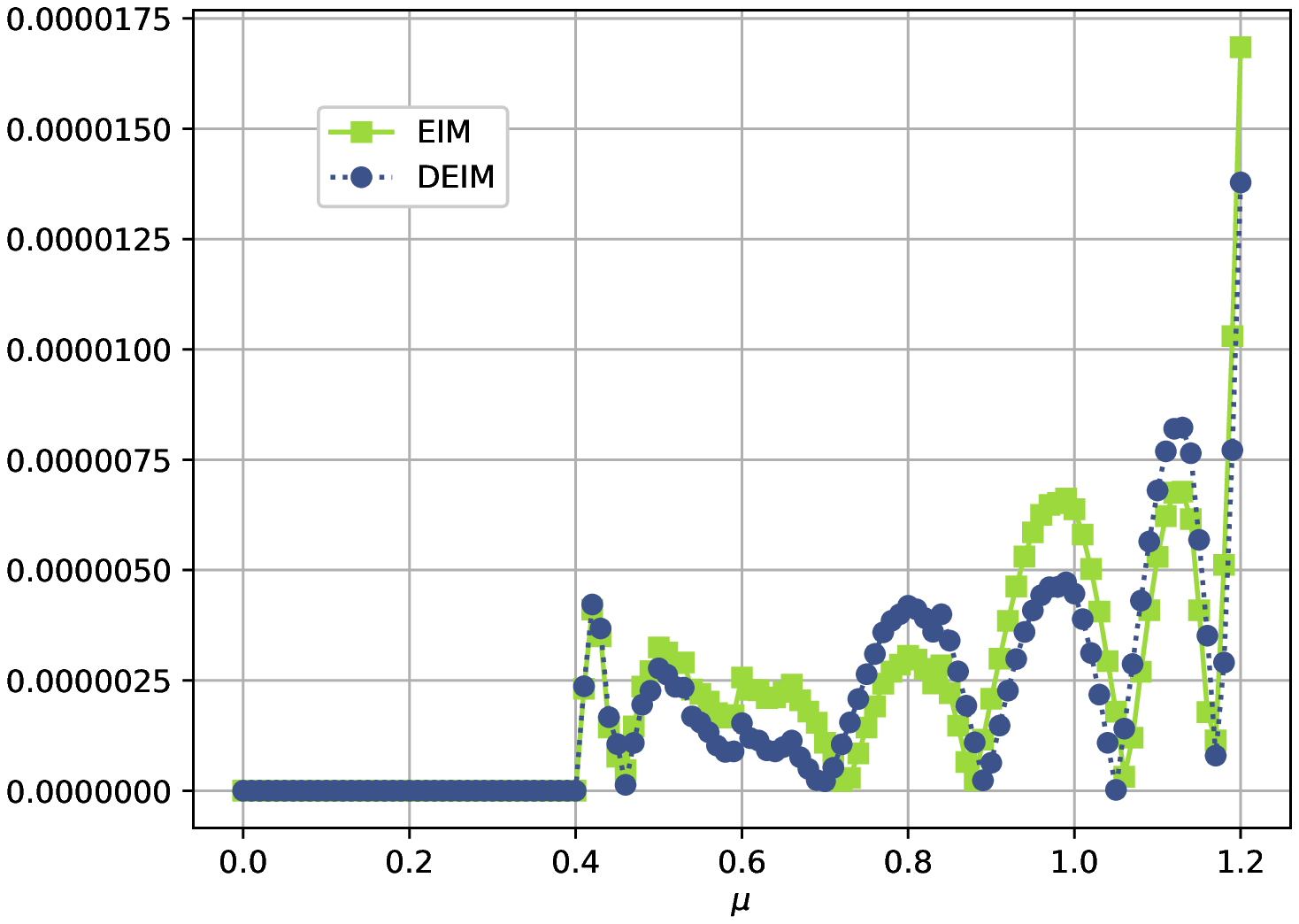}
    \put(-6,40){\scriptsize $E_N$}
  \end{overpic}
  \hfill
  \begin{overpic}[width=.48\textwidth,grid=false]{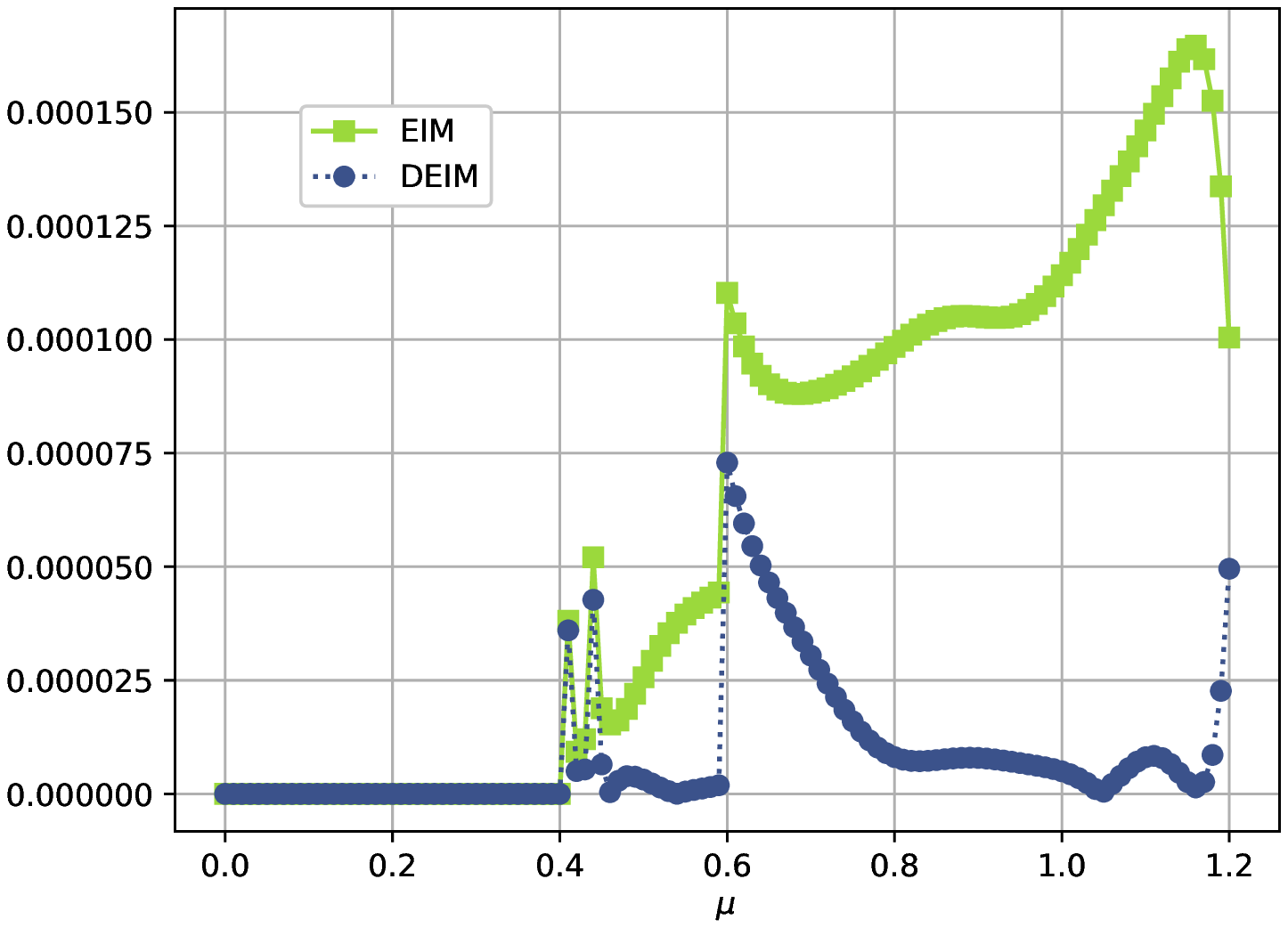}
   \put(-6,40){\scriptsize $E_\rho$}
  \end{overpic}
\caption{Difference in absolute value between the branches of the bifurcation diagram computed with FOM and ROM with EIM/DEIM in the $\mu$-$N$ plane (left) and in the $\mu$-$|| \rho||_\infty$ plane (right), i.e.~reduced order errors $E_N$ (right) and $E_\rho$ (right) for branch $|0,1\rangle$.
}\label{fig:EIM_DEIM}
\end{center}
\end{figure}

The purpose of this section is two-fold: show that the branch is reconstructed accurately and report on the substantial computational time savings enable by the hyper-reduction strategies. Fig.~\ref{fig:EIM_DEIM} shows reduced order errors $E_N$ (right) and $E_\rho$ (right) for branch $|0,1\rangle$.
We see that the values of $E_N$ for EIM and DEIM are comparable over the entire $[0, 1.2]$ interval,
with a peak at $\mu = 1.2$ of the order of $10^{-5}$. Error $E_\rho$ is slightly larger for EIM, with a peak of the order 
of $10^{-4}$.

%
%

\begin{figure}[h]
\begin{center}
\begin{overpic}[width=.48\textwidth,grid=false]{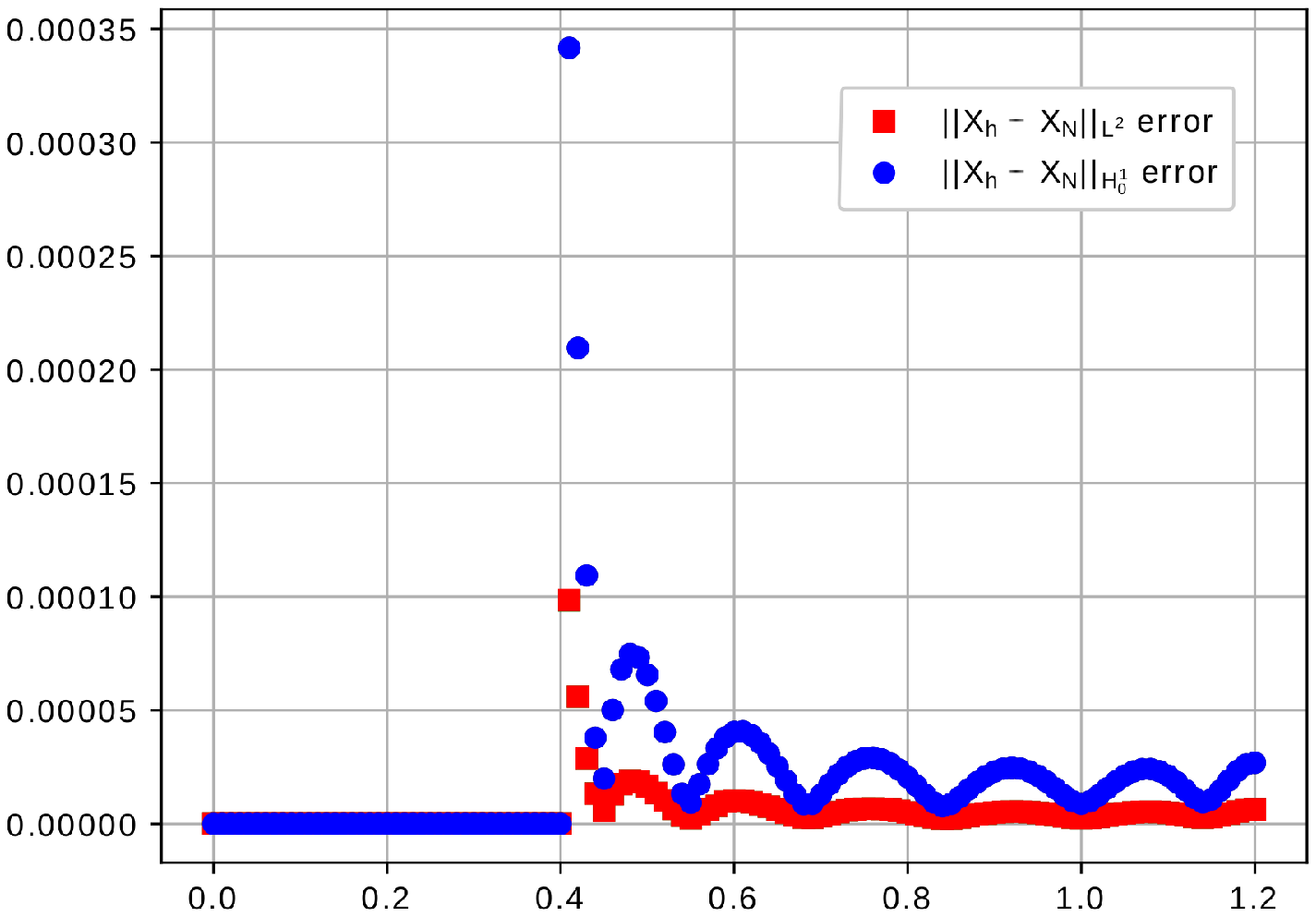}
    \put(55,2){\scriptsize $\mu$}
  \end{overpic}
  \hfill
  \begin{overpic}[width=.48\textwidth,grid=false]{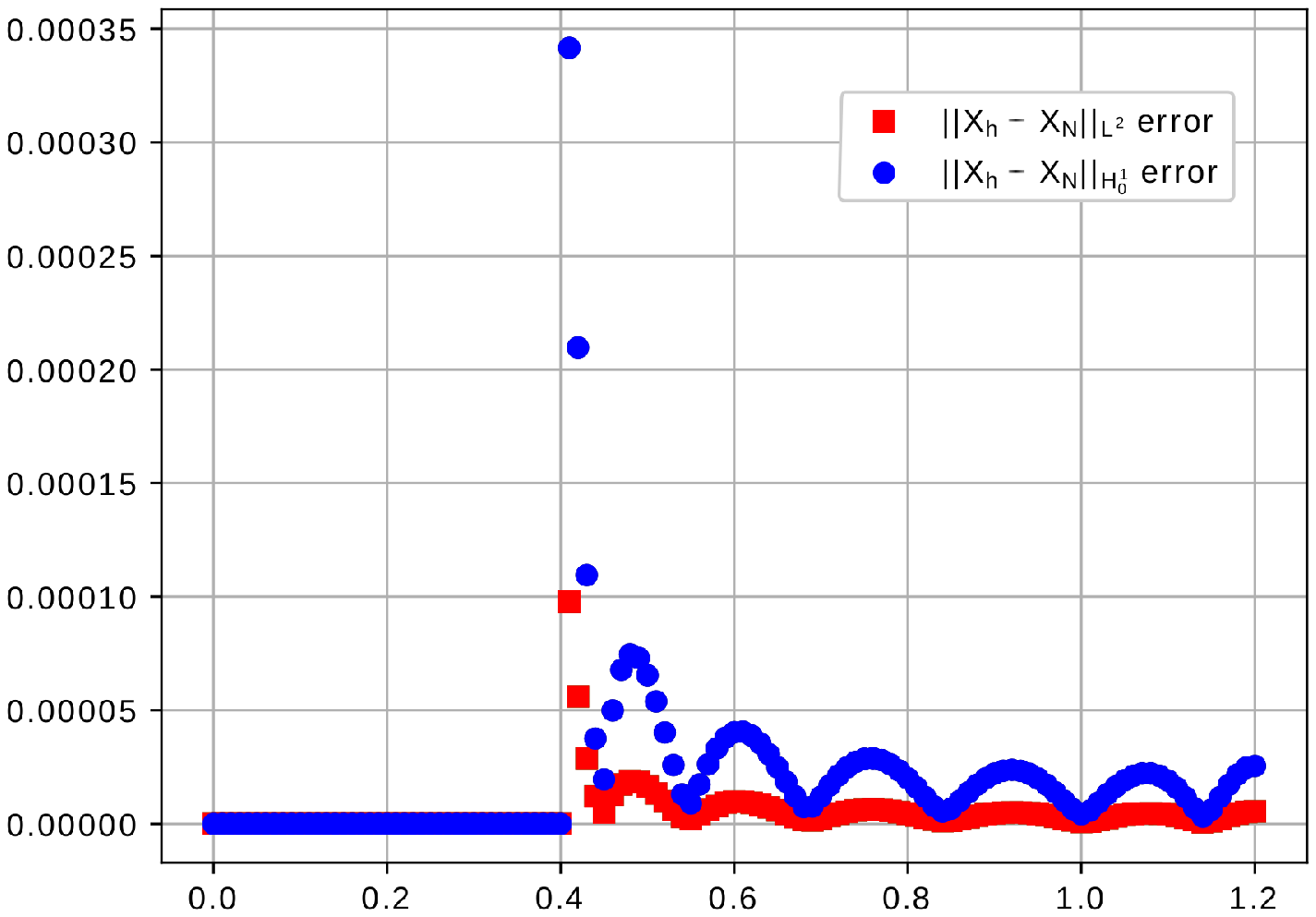}
    \put(55,2){\scriptsize $\mu$}
  \end{overpic}
\caption{Difference between $X_h$ (FOM solution) and $X_N$ (ROM solution) computed with EIM (left) and DEIM (right)
in the $L^2$ and $H^1_0$ norms for branch $|0,1\rangle$.
}\label{fig:EIM_DEIM_error}
\end{center}
\end{figure}

Finally, Fig.~\ref{fig:EIM_DEIM_error} shows the difference between FOM solution and ROM solution 
computed with EIM (left) and DEIM (right) in the $L^2$ and $H^1_0$ norms, again for branch $|0,1\rangle$.
We observe slightly larger error peaks than in the case of ROM with no affine-recovery technique:
compare Fig.~\ref{fig:EIM_DEIM_error} with Fig.~\ref{fig:rho_branch2a_err}. These increased errors
are the price to pay for a considerable computational speed-up. 
With EIM, it takes only 55 s to construct branch $|0,1\rangle$ while
it takes 246 s with FOM. So, our ROM approach coupled with EIM 
is almost five times faster than the FOM. As for DEIM, the computational time 
savings are even better: it takes only 7 s for the construction of branch $|0,1\rangle$, with corresponds to
a speed-up factor of almost 32. This drastic reduction of the computational time allowed by 
the DEIM is expected \cite{Chaturantabut2010}.


The hyper-reduction techniques is especially effective in reducing the computational time
needed for the two-parameter study presented in Sec.~\ref{sec:two_parameter}.
Our ROM approach with DEIM takes only 6 minutes to reconstruct the graph in Fig.~\ref{fig:bif_3d_FOM_1}.
This is a speed-up factor of 40 with respect to our ROM approach without DEIM
and a factor 60 with respect to the FOM. 

\section{Conclusions and perspectives}\label{sec:concl}
We have exploited reduced order methods to drastically reduce the computational time
required to trace a bifurcation diagram. 
We proposed a combination of different techniques to overcome the curse of dimensionality, mostly when more 
than one physical parameter varies.
In particular we used the standard Finite Element method as an high fidelity solver for the offline phase 
and a Proper orthogonal decomposition as a reduction technique to construct a basis for the approximation manifold.
Moreover, we implemented a simple continuation method to prevent the divergence of the Newton-Kantorovich method.  

To demonstrate the effectiveness of our approach, we studied bifurcating phenomena 
in quantum mechanics described by the Gross-Pitaevskii equation.
We were able to trace the full bifurcation diagram with high accuracy 
even when dealing with a multi-parameter context. We also highlighted the need for
an affine recovery technique (EIM/DEIM) in order to obtain important computational 
time savings. 

There are several ways in which this work could be expanded. 
One can consider different reduction strategies, such as POD-Greedy for the 2-parameters test case, in order 
to understand how the bifurcating parameter varies and then discover new branches at a reduced
computational cost.
In addition, a deflation method could be implemented and paired with a more involved and smart continuation method.
Finally, an extension to the 3D version of the same model or to multi-component systems could be investigated.

\section*{Acknowledgements}  

This work was supported by European Union Funding for Research and Innovation through the European Research Council
(project H2020 ERC CoG 2015 AROMA-CFD project 681447, P.I. Prof. G.~Rozza).
This work was also partially supported by NSF through grant DMS-1620384 (Prof. A.~Quaini).
We acknowledge fruitful conversations with Dr.~S. Garner. We thank Dr.~F.~Ballarin (SISSA) for his great help with the RBniCS software library and precious discussion.

\bibliographystyle{plain}
\bibliography{rbsissa}

\end{document}